                       \documentstyle[12pt]{amsart}
                        \numberwithin{equation}{section}
                        \theoremstyle{plain}
                                
\newtheorem{thm}{Theorem}[section]

 \newtheorem{lem}[thm]{Lemma}

\newtheorem{pro}[thm]{Proposition}

                        \theoremstyle{definition}
                                \newtheorem{rem}{Remark}[section]

\newtheorem{defn}{Definition}[section]

\input amssym.def
                
\begin{document}

    \title[A universal invariant of homology 3-spheres]
{ An invariant of integral homology 3-spheres 
\protect \\
which is universal for all finite type invariants}

                \author[ Thang  T. Q. Le]{Thang  T. Q. Le }
\address{Dept. of Mathematics, 106 Diefendorf Hall, 
 SUNY at Buffalo, Buffalo, NY 14214, USA }
\email{letu@@newton.math.buffalo.edu}

\maketitle
\begin{abstract} In [LMO] a 3-manifold invariant $\Omega(M)$ is 
constructed
using a modification of the Kontsevich integral and the Kirby 
calculus. The invariant $\Omega$ takes values in a graded Hopf
algebra of Feynman 3-valent graphs.
Here we show that for homology 3-spheres the invariant $\Omega$ 
is {\em universal} for all finite type invariants, i.e. $\Omega_n$ is 
an invariant of order $3n$ which
 dominates all other invariants of  the same order. 
This shows that the set of finite type invariants of homology
3-spheres is equivalent to the Hopf algebra of Feynman 3-valent
graphs.
Some corollaries are discussed. A theory of groups of homology
3-spheres, similar to Gusarov's theory for knots, is presented.
\end{abstract}

\newcommand{\bn}[2]{\text{$\left(\begin{array}{c}#1\\ #2\end{array}
\right)$}}

\newcommand{\cS}{{\cal S}}
\newcommand{\cA}{{\cal A}}

\newcommand{\VM}{{\cal {VM}}}
\newcommand{\cF}{{\cal F}}
\newcommand{\cP}{{\cal P}}
\newcommand{\cC}{{\cal A}}
\newcommand{\cE}{{\cal E}}
\newcommand{\hcC}{\hat{\cal A}}
\newcommand{\cZ}{{\check Z}}
\newcommand{\hZ}{{\hat  Z}}
\newcommand{\cCO}{{\cal D}}
\newcommand{\coC}{\overset{\circ}{\cal A}}
\newcommand{\Sl}{{\sqcup^lS^1}}
\newcommand{\Gr}{\text{\rm Grad}}
\newcommand{\Grn}{\text{\rm Grad}_n}
\newcommand{\bR}{{\Bbb R}}
\newcommand{\cM}{{\cal M}}
\newcommand{\On}{O^*_{3n}}

\section{Introduction}

In a previous work \cite{LMO} (joint with J. Murakami and T. Ohtsuki)
  we defined a 3-manifold invariant $\Omega$ with values in,
as expected in perturbative Chern-Simon theory (see 
\cite{AS1,AS2,Ko2,Roz,Witten}),
the graded Hopf algebra $\cCO$ of vertex-oriented 3-valent graphs. 
The construction 
uses a  modification of the Kontsevich integral and the 
Kirby calculus.
Actually, $\Omega$ corresponds only to the trivial connection 
in the
perturbation theory.  One can expect that $\Omega$ plays an
 important role in the set of homology 3-sphere invariants. 

The first non-trivial degree part $\Omega_1$ is essentially
the Casson invariant; and $\Omega$ can be regarded as a far 
generalization of
the Casson invariant. The degree $n$ part of $\Omega$ is constructed
 using
finite type invariants of links of order $\le (l+1)n$, where $l$ is
 the number
of components of the link.

Finite type invariants for {\it knots and links} were introduced by 
Vassiliev, see \cite{Vas,BL}, and proved to be very useful in knot 
theory.
Kontsevich \cite{Ko1} introduced a knot invariant, called the
 {\it Kontsevich integral}, 
with values
in a graded algebra of chord diagrams. The grading $n$ part of 
the Kontsevich integral
is a finite type invariant of order $n$; and it is universal, 
since it dominates all
other invariants of the same order.
All  quantum invariants, including the Jones polynomial, are 
special values of the
Kontsevich integral, obtained  using weight 
systems coming from
semi-simple Lie algebras (see Theorem 10 of \cite{LM2} and 
Theorem XX.8.3 of \cite{Kas}).

Using $sl_2$ quantum invariants at roots of unity, Reshetikhin 
and Turaev \cite{RT2} constructed quantum invariants of 3-manifolds.
 There are generalizations for other quantum groups.
The algebraic part of the construction is the complicated theory 
of {\it modular}
Hopf algebra and quantum groups {\it at roots of unity}, see 
\cite{Tur}.

Since quantum invariants of links are special values of the 
Kontsevich integral,
there arose a question of constructing 3-manifold invariants  
from the 
Kontsevich integral. In \cite{LMO}, we succeeded in constructing 
such an invariant $\Omega$.
Actually we used a modification of the Kontsevich integral, called
 the universal Kontsevich-Vassiliev invariant
for {\it framed} links and tangles (see \cite{LM1,LM2}).

Invariants of finite type for {\it homology 3-spheres} were 
introduced by Ohtsuki \cite{Ohtsuki},
in analogy with the knot case. The theory was developed further 
in \cite{GL,GL2,GO,GO2,Hab,Lin} and others.
In particular, it was proved that orders of finite type invariants 
are multiples of 3
(see \cite{GL,GO}), and that the space of invariants of a fixed 
order is finite-dimensional (see \cite{Ohtsuki}). 

Both finite type invariants for links and for homology 3-spheres 
can be considered
as a way to organize the set of invariants in a systematical way. 
While a lot
has been known about the knot case, much less is for the manifold 
case. 

The aim of this paper is to show that for homology 3-spheres, 
$\Omega$ is
as powerful as the set of all finite type invariants; i.e. 
$\Omega$ is
a counterpart of the Kontsevich integral in knot theory.  
We show that the 
degree $n$ part of $\Omega$ is an invariant of degree $3n$ 
which dominates 
all other invariants of degree $3n$. This implies that the set of
finite type invariants of homology 3-spheres is  equivalent
to the (purely combinatorial) Hopf algebra $\cCO$ of Feynman 3-valent 
diagrams. For example,
the number of linearly independent invariants of degree $\le n$ is
equal the number of linearly independent 3-valent diagrams of degree
$\le n$. The study of finite type invariants of homology 3-spheres
is reduced to the study of $\Omega$.

The proof of the main theorem is much more complicated than that of 
the similar theorem in the knot case. 
We have to
use many results of the theory of the (modified) Kontsevich integral. 
We  show that
every weight system, defined in \cite{GO}, can be integrated to a 
finite type invariant, like in the knot case (and hence give an 
affirmative answer to Question 1 there). 

We also show that for every $n$, there is a homology 3-sphere $M$ 
which cannot be distinguished from $S^3$ by any invariant of order 
$< 3n$, but can be distinguished
from $S^3$ by an invariant of order $3n$. An operation on homology
3-spheres which does not alter values of invariants of order less 
than a fixed number
is presented.
 A theory
of {\it groups of homology 3-spheres}, similar to Gusarov's groups
of knots, is sketched.
Some other corollaries of the main theorem are presented.

The paper is organized as follows. In Sections \ref{prep}
 we recall the universal Kontsevich-Vassiliev
invariant for {\it framed oriented links and tangles}, considered in 
\cite{LM1,LM2}. We also show
some properties of this invariant, needed later. In Section 3 we 
give a definition of
$\Omega$, and recall some fundamental concepts in the theory of 
finite type invariants of
homology 3-spheres. In Section 4 we present the main results and 
discuss some corollaries and related problems. In Section 5 we prove 
the lemmas needed in the proof of the main theorem.

The author would like to thank W. Menasco, J. Murakami and T. Ohtsuki
 for useful discussions.

\section{Chord diagrams and finite type invariants for framed tangles}
\label{prep}
\subsection{Chord diagrams} Note that our definition of chord 
diagrams 
is more general than that of \cite{BN1,LM1,LM2}. All vector spaces
are over the field $\Bbb Q$ of rational numbers.

A {\it uni-trivalent graph}
is a graph every vertex of which is either univalent
or trivalent. A uni-trivalent graph is {\it vertex-oriented} if at
 each trivalent
 vertex a cyclic order of edges is fixed. A 3-valent (resp. 1-valent)
vertex is called  an {\it internal} ({\it external}) vertex.

Let $X$ be a compact oriented 1-dimensional manifold
whose components are {\it numbered}.
A {\it chord diagram} with support $X$ is the manifold $X$ together
 with a
vertex-oriented uni-trivalent graph
 whose external vertices are on $X$; and the graph does not have
 any connected component homeomorphic to a circle.

In  figures components of $X$ are depicted by solid lines,
while the graph is depicted by dashed lines, with the convention that
the orientation at every vertex is counterclockwise.
Each chord diagram has a natural topology.
Two chord diagrams $\xi,\xi'$ on $X$ are regarded
as equal if there is a homeomorphism
$f:\xi\to \xi'$ such that $f|_{X}$ is a homeomorphism
of $X$ which preserves components and orientation and the restriction 
of
$f$ to the dashed graph is a homeomorphism preserving orientation at 
every vertex.

 There may be connected components of the  
dashed graph which do
not have univalent vertices, and hence do not connect to any solid 
lines.

Let $\cC(X)$ be the vector space 
spanned by chord diagrams with support $X$,
subject to the AS, IHX and STU relations.
The AS ({\it anti-symmetry}) condition says that $\xi_1+\xi_2=0$ for 
any two chord diagrams which are identical except for the orientation
at a vertex. The IHX (or Jacobi) relation says that 
 $\xi_1=\xi_2+\xi_3$, for every three chord diagrams identical 
outside a ball in which they differ as in Figure \ref{IHX}.

\begin{figure}[htp]
\centerline{
\begin{picture}(60,60)
\multiput(0,20)(3,0){20}{\line(1,0){2}}
\multiput(0,40)(3,0){5}{\line(1,0){2}}
\multiput(15,40)(0,-3){7}{\line(0,-1){2}}
\multiput(30,60)(0,-3){14}{\line(0,-1){2}}
\put(30,3){\makebox(0,0){$\xi_1$}}
\put(15,20){\circle*{2.5}}
\put(30,20){\circle*{2.5}}
\end{picture}\hskip 1.5cm
\begin{picture}(60,60)
\multiput(0,20)(3,0){20}{\line(1,0){2}}
\multiput(0,40)(3,0){10}{\line(1,0){2}}
\multiput(30,60)(0,-3){14}{\line(0,-1){2}}
\put(30,3){\makebox(0,0){$\xi_2$}}
\put(30,20){\circle*{2.5}}
\put(30,40){\circle*{2.5}}
\end{picture}\hskip 1.5cm
\begin{picture}(60,60)
\multiput(0,20)(3,0){20}{\line(1,0){2}}
\multiput(0,40)(3,0){15}{\line(1,0){2}}
\multiput(45,40)(0,-3){7}{\line(0,-1){2}}
\multiput(30,60)(0,-3){14}{\line(0,-1){2}}
\put(30,3){\makebox(0,0){$\xi_3$}}
\put(30,20){\circle*{2.5}}
\put(45,20){\circle*{2.5}}
\end{picture}}
\caption{\label{IHX}}\end{figure}

\begin{figure}[htp]

\centerline{
\begin{picture}(60,40)
\multiput(20,40)(0,-3){7}{\line(0,1){2}}
\multiput(40,40)(0,-3){7}{\line(0,1){2}}
\multiput(30,18.5)(0,-3){7}{\line(0,1){2}}
\multiput(20,20)(3,0){7}{\line(1,0){2}}
\put(30,20){\circle*{2.5}}
\thicklines
\put(0,0){\vector(1,0){60}}
\end{picture}\hskip 1cm
$=$\hskip 1cm
\begin{picture}(60,40)
\multiput(20,40)(0,-3){14}{\line(0,1){2}}
\multiput(40,40)(0,-3){14}{\line(0,1){2}}
\thicklines
\put(0,0){\vector(1,0){60}}
\end{picture}\hskip 0.6cm
$-$\hskip 0.6cm
\begin{picture}(60,40)
\multiput(20,40)(0,-3){7}{\line(0,-1){2}}
\multiput(40,20)(0,-3){7}{\line(0,-1){2}}
\multiput(20,-.5)(0,3){4}{\line(0,1){2}}
\multiput(30,10)(0,3){7}{\line(0,1){2}}
\multiput(40,30)(0,3){4}{\line(0,1){2}}
\multiput(19.5,10)(3,0){4}{\line(1,0){2}}
\multiput(29.5,30)(3,0){4}{\line(1,0){2}}
\multiput(20,20)(3,0){3}{\line(1,0){2}}
\multiput(32,20)(3,0){3}{\line(1,0){2}}
\thicklines
\put(0,0){\vector(1,0){60}}
\end{picture}}

\caption{The STU relation\label{STU}}\end{figure}
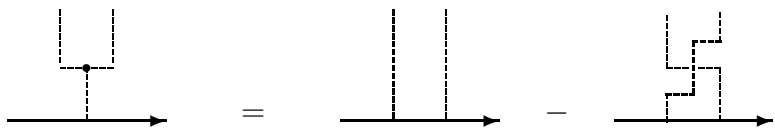

The STU relation is depicted in Figures \ref{STU}.

The {\it degree} of a chord diagram is
half the number of (external and internal)  vertices of the {\it 
dashed graph}.
Since the relations AS, IHX and STU respect the degree,
there is a grading on $\cC(X)$ induced by this degree.
We also use  $\cC(X)$ to denote the completion of $\cC(X)$
with respect to the grading.

Let $\cCO_n$ be the vector space spanned by all vertex-oriented 
3-valent 
(without 1-valent vertex) graphs
of degree $n$ 
subject to the AS and IHX relations. Let $\cCO_0$ be the ground field
 $\Bbb Q$, and
let $\cCO$ be the direct product of all $\cCO_n, n=0,1,2,\dots$.
This vector space $\cCO$ is
a commutative algebra in which the product of 2 graphs is just the
 union
 of them.
All $\cC(X)$ can be regarded as graded $\cCO$-modules, where the 
product of 
a 3-valent graph $\xi\in\cCO$ and a chord diagram $\xi'\in\cC(X)$ 
is their 
disjoint union. In \cite{LMO}, $\cCO$ is denoted by $\cA(\emptyset)$.

Suppose $C$ is a component of $X$.
Reversing the orientation of $C$, from $X$  we
get $X'$. Let $S_{(C)}: \cA(X)\to\cA(X')$ be the linear mapping which
transfers every chord diagram $\xi$ in $\cA(X)$ to $S_{(C)}(\xi)$ 
obtained
from $\xi$ by reversing the orientation of $C$ and multiplying  by
$(-1)^m$, where $m$ is the number of vertices of the dashed graph
on the component $C$.

We also define the linear operator $\Delta_{(C)}$ which is described
in Figure \ref{Delta}. The definition is explained as follows.
\begin{figure}[htp]
\centerline{
\begin{picture}(100,60)
\put(-15,40){\dashbox{2}(15,0){}}
\thicklines \put(0,60){\line(0,-1){40}}
\put(0,5){\makebox(0,0){$C$}}
\put(50,40){\makebox(0,0){$\to$}}
\end{picture}
\begin{picture}(120,60)
\put(-15,40){\dashbox{2}(15,0){}}
\put(100,40){\dashbox{2}(20,0){}}
\put(50,40){\makebox(0,0){$+$}}
\thicklines \put(0,60){\line(0,-1){40}}
\put(5,20){\line(0,1){40}}
\put(120,60){\line(0,-1){40}}
\put(115,20){\line(0,1){40}}
\end{picture}
}
\caption{\label{Delta}}\end{figure}

Replacing $C$ by $2$ copies of $C$,
from $X$ we get $X^{(2,C)}$,
with a projection $p:X^{(2,C)}\to X$.
If $x$ is a point on $C$ then $p^{-1}(x)$ consists of $2$ points,
while if $x$ is a point of other components,
then $p^{-1}(x)$ consists of one point.
Let $\xi$ be a chord diagram on $X$,
with the dashed graph $G$.
Suppose that there are $m$ univalent vertices of $G$ on $C$.
Consider all possible new chord diagrams on $X^{(2,C)}$
with the same dashed graph $G$ such that
if a univalent vertex of $G$ is attached to a point $x$ on $X$ in 
$\xi$,
then this vertex is attached to a point in $p^{-1}(x)$
in the new chord diagram.
There are $2^m$ such chord diagrams,
and their sum is denoted by $\Delta_{(C)}(\xi)$.
It is easy to check that these linear mappings $S_{(C)}$ 
and $\Delta_{(C)}$ are well-defined.

Suppose that $X$ has a distinguished component
$C$.
Let $\xi\in\cA(X)$ and $\xi'\in\cA(S^1)$ be two chord diagrams.
{}From each  of $C$ and $S^1$ we remove a small arc
which does not contain any vertices.
The remaining part of $S^1$ is an arc which we glue  to $C$ in
the place of the removed arc so that the orientations are 
compatible.
The new chord diagram is called the
{\it connected sum of $\xi$ and $\xi'$ along $C$};
it does not depend on the
locations of the removed arcs. The proof is
the same as in the case $X=S^1$, considered in \cite{BN1}.

When $X=S^1$, the connected sum defines an algebra structure 
on $\cA(S^1)$ which is known to be commutative (see \cite{BN1}).

We now define a co-multiplication $\hat \Delta$
in $\cC(X)$.
A {\it  chord sub-diagram} of a chord diagram $\xi$ with dashed
 graph $G$
 is any chord
 diagram obtained from $\xi$ by removing
some  connected components of $G$. The {\it complement chord 
sub-diagram}
of  a chord sub-diagram
$\xi'$
is  the chord  sub-diagram obtained by removing components of $G$
which are in $\xi'$.
Let
$$\hat \Delta (\xi)=\sum \xi'\otimes \xi''.$$
Here the sum is over all chord sub-diagrams $\xi'$ of $\xi$, 
and $\xi''$ is
 the complement
of $\xi'$. This co-multiplication is co-commutative.

A similar co-product is defined for $\cCO$, by the same formula,
 and 
$\cCO$ becomes
a commutative co-commutative Hopf algebra. Hence $\cCO$ is the
 polynomial algebra
on {\it primitive} elements, i.e. elements $x$ such that 
${\hat\Delta}(x)=1\otimes x + x\otimes 1$. It is easy to 
see that an element is primitive if and only
if it is a linear combination of {\it connected} 3-valent graphs.

\subsection{A universal invariant for framed tangles: the modified 
Kontsevich integral}

We  fix an oriented 3-dimensional
 Euclidean space $\bR^3$ with coordinates $(x,y,t)$.
{\it A tangle} is a smooth one-dimensional compact  oriented 
manifold
 $L\subset \bR^3$  lying between two horizontal planes 
$\{t=a\},\{t=b\},
 a<b$ such that all the boundary points are lying on two lines 
$\{t=a,y=0\},
\{t=b,y=0\}$, and at every boundary point $L$ is orthogonal to these 
two planes. These lines are called the top and the bottom lines
 of the tangle.

{\it A normal vector field} on a tangle $L$ is a smooth vector field 
on
$L$ which is nowhere tangent to $L$ (and, in particular, is nowhere
zero) and which is given by the vector $(0,-1,0)$ at every boundary
point. {\it A framed tangle} is a tangle enhanced with a normal vector
field. Two framed tangles are isotopic if they can be deformed by  a
1-parameter family of diffeomorphisms into one another within the 
class
of framed tangles.

Framed oriented links are special framed tangles when there is no 
boundary point.
The empty link, or empty tangle, by definition, is the empty set.

One  assigns a symbol $+$ or $-$ to all the boundary points of a
tangle according to  whether the tangent vector at this point directs
downwards or upwards. Then on the top boundary line of a tangle 
diagram
we have a word $w_t$ of symbols consisting of $+$ and $-$. Similarly 
on the bottom
boundary line there is also a word $w_b$ of symbols $+$ and $-$.

A non-associative word on $+,-$ is an element of the {\it free magma}
generated by $+,-$ (see \cite{Ser}). A {\it $q$-tangle} (or 
{\it non-associative}
tangle) is a tangle whose top and bottom words $w_t,w_b$ are 
equipped with
some non-associative structure.

Every tangle in this paper will also be regarded as a $q$-tangle, 
with {\it standard} non-associative structure on $w_t,w_b$: 
 product is taking successively from left to right.

If $T_1, T_2$ are tangles such that $w_b(T_1)=w_t(T_2)$
we can define the product $T=T_1\,T_2$ by placing $T_1$ on 
top of
$T_2$. In this case,
if  $\xi_1\in 
\cA(T_1),
 \xi_2\in \cA(T_2)$ are chord diagrams, then the {\it product} 
$\xi_1\xi_2$
 is a chord diagram in $\cA(T)$ obtained by placing $\xi_1$ on 
top of $\xi_2$.

For any two tangles $T_1,T_2$ with the same top and bottom lines, 
we can define
their {\it tensor product} $T_1\otimes  T_2$ by putting $T_2$ to 
the right
 of $T_1$.
Similarly, if $\xi_1\in\cC(T_1),\xi_2\in\cC(T_2)$ are chord diagrams, 
then one defines $\xi_1\otimes \xi_2\in \cC(T_1\otimes  T_2)$ 
by the same way.

For example, if $L_1, L_2$ are two links, then $L_1L_2=L_1\otimes 
L_2=L_1\sqcup L_2$,
 the disjoint union of $L_1$ and $L_2$.

Let $T$ be a framed tangle, then there exists an element 
$\hZ(T)\in\cC(T)$ 
called the {\it universal Kontsevich-Vassiliev invariant} of $T$. 
This 
invariant satisfies 
\begin{equation}
\text{if $T=T_1T_2$ then $\hZ(T)=\hZ(T_1)\hZ(T_2)$}.
\label{multi}
\end{equation}

The definition of this invariant $\hZ$ is given in, say,  
\cite{LM1,LM2} (see also \cite{BN2}); 
in \cite{LM1,LM2} $\hZ$ is denoted by 
$\hZ_f$. It is a regularization of the original Kontsevich integral,
regularized so that it can be defined for {\it framed links}. For 
framed links,
 $\hZ$ is
equivalent to the set of all finite type invariants.

Note that the definition of $\hZ$ depends on an element called {\it
associator}. For the theory of associators in our sense, see, for 
example, \cite{LM2}.
 For links the values of $\hZ$ do not depend on what associator to 
choose 
(see Theorems 7,8 in \cite{LM2}), but the values of $\hZ$  of a 
framed tangle with non-empty boundary
do depend on the associator. For technical convenience,
we will use the
normalization of $\hZ$ as in \cite{LM3} and 
the associator in \cite{LM3}, since 
this associator has rational coefficients and admits a lot a symmetry 
(more symmetry than the 
associator coming from the Knizhnik-Zamolodchikov equations).
Note that the normalization of $\hZ$ is chosen so that $\hZ$ of an 
empty link is 1.

Let $\cP_n$ be the space of chord diagrams with support being $n$
 numbered straight vertical
lines, pointing downwards. There is a structure of a (non-commutative,
if $n>1$,) 
algebra on $\cP_n$,
where the product of two chord diagrams is obtained by placing the 
first on top of the second.
Together with the co-multiplication $\hat\Delta$, we get a 
co-commutative Hopf algebra $\cP_n$.   It is known that
the subspace of  
primitive elements in $\cP_n$ is spanned by chord diagrams with 
{\it connected 
dashed graph}.

In what follows, element $1\in\cC(X)$ is always the chord diagram
without dashed graph.
\begin{defn}{\it 
We say that an element $\xi\in\cC(X)$ has {\it {\rm i}-filter}
 $\ge n$ if $\xi$ 
is a linear combination of
chord diagrams, each has at least $n$ internal vertices.}
\end{defn}
If $\xi$ is a chord diagram with {\it connected} dashed graph and
of degree $n$, then it is easy to see that $\xi$ has at least
$(n-1)$ internal vertices. Hence we have the following
\begin{pro}
If $\xi\in \cP_n$ is primitive and is a linear combination of chord
diagrams of degree $\ge n$, then $\xi$ has {\rm i}-filter $\ge (n-1)$.
\label{ifilter}
\end{pro}
It is known that the associator is an element in $\cP_3$
of the form $\exp(\xi)$, where $\xi$ is {\it
primitive} and is a linear combination of chord diagrams of
degree $\ge 2$. It follows that
\begin{lem} The associator has the form $1+(\text{elements of 
 {\rm i}-filter 
$\ge 1$})$.
\end{lem}
Unlike invariants of tangles coming from co-associative quantum 
groups,
in general, if  $T=T_1\otimes T_2$ then $\hZ(T)\not=\hZ(T_1)
\otimes\hZ(T_2)$.
 But one has
\begin{lem} For some elements 
$a,b\in\cC(T_1\otimes T_2)$  of the form $1+\text{(elements of  {\rm i}-filter 
$\ge 1$)}$, one has
$$\hZ(T_1\otimes T_2)=a[\hZ(T_1)\otimes\hZ(T_2)]b.$$
 Moreover, if $T_2$ is a trivial tangle, then $a=b=1$.
\label{hZ}
\end{lem}
\begin{pf} The general method of computing $\hZ$ in \cite{LM2} 
says that
$$\hZ(T_1\otimes T_2 )=a[\hZ(T_1)\otimes\hZ(T_2)]b,$$ where $a,b$
 are derived from the associator
by some operators; and these operators do not decrease 
{\rm i}-filters.
\end{pf}

For a framed tangle $T$ and its component $C$, let
  $T^{(2,C)}$ be obtained from $T$ by adding a string
 $C^\prime$ parallel to $C$ with respect to the framing.

\begin{thm}\label{theorem1}\cite{LM3}
Let $C$ be a component of a framed tangle $T$.

a) One has $$\hat Z(T^{(2,C)})=\Delta_{(C)}(\hat Z(T)) + 
\text{(elements of i-filter $\ge 1$)}.$$ 

b) Let $T'$ be obtained from $T$ by reversing the orientation of 
$C$, then
$$\hat Z(T')=S_{(C)}(\hat Z(T)).$$
\end{thm}
Actually,  Theorem 4.2 of \cite{LM3} (the proof of which requires 
the special 
associator mentioned above) says that we have the identity in $a)$, 
without
any extra terms of  elements of i-filter $\ge 1$, but for another 
non-associative
structure  of $T^{(2,C)}$. Here we have the standard order; and 
the values of $\hZ$
of two tangles differed by non-associative structure 
differ by an element of i-filter $\ge 1$, by the same reason as in 
the proof of Lemma \ref{hZ}. Part $b)$
follows easily from the definition of $\hZ$, see Theorem 4 of 
\cite{LM2}.

\section{The universal invariant $\Omega(M)$; finite type invariants}

For a graded space $A$, we will denote by $\Grn(A)$ the subspace of
 grading 
$n$; by $\Gr_{\le n}(A)$ the subspace of grading $\le n$. For an 
element $x\in A$,
 we will denote by $\Grn(x)$ and $\Gr_{\le n}(x)$, respectively, 
the projection
 of $x$ on $\Grn(A)$ and 
$\Gr_{\le n}(A)$.

\subsection{The mapping $\iota_n$}\label{iota}
The mappings $\iota_n$ plays an important role in the theory; it helps 
to convert
chord diagrams with  supports to chord diagrams without  support.

We denote by $\coC(X)$ the space of chord diagrams with support $X$, 
subject to
 the AS, IHX and STU relations, as in $\cC(X)$; but {\it the dashed
 graph may 
contain some finite number of dashed components
which are loops}. Certainly $\cC(X)$ is a subspace of $\coC(X)$, 
and both are
graded  $\cCO$-modules.
 
Let $(L<2n)$ be the equivalence relation in $\coC(X)$ generated by: a 
chord diagram with less than $2n$ external vertices on some solid 
component
 is equivalent to 0.

Suppose on the boundary of a ball $B$ there are $2n$ distinct points. 
Connect every point
of these $2n$ points with exactly one other point by a dashed line; 
there will be
 $(2n-1)!!$ ways to
do that, each produces a 1-valent dashed graph in $B$. The formal sum 
of 
these $(2n-1)!!$ 
1-valent graphs is denoted by $T^n_{2n}$.
Let $R_n$ be the equivalence relation $T^n_{2n}=0$; i.e. if we have 
$(2n-1)!!$ 
chord diagrams identical everywhere, except for a ball $B$ in which 
they are
 exactly terms of $T^n_{2n}$, then the sum of these chord diagrams 
is 0
(in \cite{LMO}, $R_n$ is denoted by $P_n$).

Let $O_n$ be the equivalence relation in $\coC(X)$ generated by: if 
$\xi$ is a chord 
diagram which is a union of $\xi'$ and a dashed component which is 
just a
 loop, then $\xi=-2n\, \xi'$. With this relation one can remove all 
loop 
dashed components.
\begin{pro} The quotient space of
$\Gr_{\le n(l+1)}[\coC(\Sl)]$, where $\Sl$ is the union of $l$ solid 
loops,
 by relations $(L<2n)$, $R_{n+1}$ and $O_n$ is the 
free $\Gr_{\le n}(\cCO)$-module generated by $(x_n)^{\otimes l}$,
 where 
$(x_n)^{\otimes l}$ is the chord diagram which is the disjoint 
union of 
$l$ solid loops, each has $n$ isolated chords.
\end{pro}
Here an isolated chord is a dashed line (without any 3-valent vertex 
on it) 
connecting two neighboring external vertices on a solid line or loop. 
This proposition was proved in \cite{LMO} (see Lemma 3.1 there).

Suppose  $\xi$ is a chord diagram in  $\coC(\Sl)$.
 Projecting $\xi$ to
$\Gr_{\le n(l+1)}[\coC(\Sl)]/(L<2n)(R_{n+1})(O_n)$, we get an element
 $\tilde{\iota}_n(\xi) \times (x_n)^{\otimes l}$.
Here $\tilde{\iota}_n(\xi)\in \Gr_{\le n}(\cCO)$. Let 
$$\iota_n(\xi)=[(-2)^nn!]^l\,\, \tilde \iota_n(\xi)
\in \Gr_{\le n}(\cCO).$$

Note that $\iota_n$ decreases the degree of chord diagrams by $ln$,
where $l$ is the number of solid loops. Another, more geometric, 
definition of $\iota_n$ is given in \S\ref{uni-inva}.

\subsection{Definition of the universal invariant $\Omega$;
 some properties}
Suppose $M$ is an oriented closed 3-manifold obtained from $S^3$ by 
surgery 
on a framed $l$-component unoriented  link $L$. Providing $L$ with an 
arbitrary 
orientation, we can define $\hZ(L)$.
Let
$$\cZ(L)=(\nu^{\otimes l})\# \hZ(L).$$
This means, $\cZ(L)$ is obtained from $\hZ(L)$ by successively  taking 
connected 
sum of $\hZ(L)$ with
$\nu$ along each component of $L$. Here $\nu\in\cC(S^1)$ is 
 $\hZ$ of the unknot with framing 0.

We will construct an invariant $\Omega(M)\in\cCO$, using $\cZ(L)$. 
The degree $n$ part $\Grn(\Omega)$ is constructed
using  $\Gr_{\le (l+1)n}[\hZ(L)]$.

Let $U_+$  (resp. $U_-$) be the unknot with framing $+1$ 
(resp. $-1$). 
Suppose the linking matrix of $L$ has $\sigma_+$ positive eigenvalues 
and 
$\sigma_-$ negative eigenvalues.
Define (\cite{LMO})
\begin{equation}
\Omega_n(L)=\frac{\iota_n(\cZ(L))}{ (\iota_n(\cZ(U_+)))^{-\sigma_+}
(\iota_n(\cZ(U_-)))
^{-\sigma_-}}\in \Gr_{\le n}(\cCO).\end{equation}
In \cite{LMO} we proved that $ \iota_n(\cZ(U_{\pm}))$ 
are of the form 
$(\mp 1)^n+(\text{elements of degree} \ge 1)$, hence their inverses 
exist.

\begin{thm}{\rm (\cite{LMO})} $\Omega_n(L)$ does not depend on the 
orientation 
of $L$ and does not change under the Kirby moves.
Hence  $\Omega_n(L)$ is an invariant of the 3-manifold $M$.
\end{thm}
Recall that $\cCO$ is a commutative co-commutative Hopf algebra.
Let
$$\Omega(M)=1+\Gr_1(\Omega_1(M))+\dots+\Gr_n(\Omega_n(M))+\dots 
\in \cCO.$$
\begin{pro}{\rm (\cite{LMO})}\label{grouplike}\quad
$\Omega(M)$ is a group-like element, i.e.
  $$\hat \Delta (\Omega (M))=\Omega(M)\otimes \Omega(M).$$ 
Hence 
 $\ln(\Omega(M))$ is a linear combination of {\em connected} 3-valent 
vertex-oriented graphs.
\end{pro}
In general, $\Omega_n(M)$ is not equal to $\Gr_{\le n}(\Omega(M))$. 
Let $d(M)$ be the
 cardinality of
$H_1(M,{\Bbb Z})$ if the first Betti number of $M$ is 0, otherwise 
let $d(M)=0$.
\begin{pro}{\rm (\cite{LMO})} We have that $\Gr_{\le n}\Omega_{n+1}(M)=d(M)\, 
\Omega_n(M).$
Hence if $M$ is an integral homology 3-sphere, then
$$\Omega_n(M)=\Gr_{\le n}(\Omega(M)).$$
\label{ax}
\end{pro}
It was proved in \cite{LMMO} that $\Omega_1$ is, in essential,
 the Casson invariant.
\begin{pro}{\rm (\cite{LMO})}\label{connectedsum}
If  $M_1,M_2$ are {\em integral homology 3-spheres}, then 
$$\Omega(M_1\# M_2)=\Omega(M_1) \times \Omega(M_2).$$
\end{pro}
In general, if $M_1,M_2$ are {\it rational} homology 3-spheres, 
then Proposition ~\ref{connectedsum} does 
not hold true. However, if we modify $\Omega$:
$$\Omega'(M)=1+\frac{\Gr_1(\Omega_1(M))}{d(M)}+\dots+
\frac{\Gr_n(\Omega_n(M))}
{d(M)^n}+\dots$$
Then $\Omega'(M_1\# M_2)=\Omega'(M_1) \times \Omega'(M_2)$, and 
$\hat \Delta (\Omega' (M))=\Omega'(M)\otimes \Omega'(M)$.

For every Lie algebra $\frak g$ with a non-degenerate {\it invariant}
bilinear form 
(say, a
semi-simple Lie algebra), one can define a weight system 
(see, for example, \cite{BN1,Ko1})
which transforms $\Omega(M)$
into a formal power series on a variable $h$. The relation of this 
formal power series
with quantum invariants coming from quantum groups at roots of unity 
is yet to discover.

\subsection{ Another definition of $\iota_n$} \label{uni-inva}

In \cite{LMO}, we also use another definition of $\iota_n$ which is 
more geometric. 
Actually, we first came to this form of the definition.
The aim of $\iota_n$
is to remove solid loops and replace them by appropriate dashed 
graphs.

Let $\cA(m)$, for any positive integer $m$, be the vector space 
spanned by uni-trivalent
vertex-oriented graphs with exactly $m$ 1-valent vertices located at
 $m$ numbered points, 
subject to the AS and IHX relations. The vector space $\cA(m)$ 
is a $\cCO$-module.

We would like to replace a solid loop with $m$ external 
vertices on it
by an element $\eta_m\in\cA(m)$, i.e. we remove the solid
 loop, then glue the
external vertices of $\eta_m$ to the 1-valent vertices of 
the removed solid loop.
Then the sequence $\eta_m$  must satisfy 

\begin{equation}
\begin{picture}(60,70)(0,-38)
\multiput(20,40)(0,-3){7}{\line(0,1){2}}
\multiput(40,40)(0,-3){7}{\line(0,1){2}}
\multiput(30,18.5)(0,-3){7}{\line(0,1){2}}
\multiput(20,20)(3,0){7}{\line(1,0){2}}
\put(30,20){\circle*{2.5}}
\put(30,-20){\circle{40}}
\put(30,-18){\makebox(0,0){$\eta_{m-1}$}}
\end{picture}\quad\quad
=\quad\quad
\begin{picture}(60,70)(0,-38)
\multiput(20,40)(0,-3){14}{\line(0,1){2}}
\multiput(40,40)(0,-3){14}{\line(0,1){2}}
\put(30,-17){\circle{40}}
\put(30,-18){\makebox(0,0){$\eta_m$}}
\end{picture}\quad - \quad
\begin{picture}(60,70)(0,-38)
\multiput(20,40)(0,-3){7}{\line(0,-1){2}}
\multiput(40,20)(0,-3){7}{\line(0,-1){2}}
\multiput(20,-.5)(0,3){4}{\line(0,1){2}}
\multiput(30,10)(0,3){7}{\line(0,1){2}}
\multiput(40,30)(0,3){4}{\line(0,1){2}}
\multiput(19.5,10)(3,0){4}{\line(1,0){2}}
\multiput(29.5,30)(3,0){4}{\line(1,0){2}}
\multiput(20,20)(3,0){3}{\line(1,0){2}}
\multiput(32,20)(3,0){3}{\line(1,0){2}}
\put(30,-18){\circle{40}}
\put(30,-18){\makebox(0,0){$\eta_m$}}
\end{picture}
\tag{p1}\end{equation}
which follows from the STU relation (we do not draw the other 
$m-2$ external vertices
of $\eta_m,\eta_{m-1}$).
In addition, 
\begin{equation}\text{$\eta_m$ is 
invariant under cyclic permutation of the external vertices}\tag{p2}.
\end{equation}

But it seems impossible to have a sequence $\eta_m\in\cA(m),m=1,2,\dots$,
 with the
help of which one can directly transform $\hZ$ into something which 
is invariant under
the Kirby moves. What we did in \cite{LMO} is, motivated by low degree 
cases and other
thoughts, to construct the transformation from $\cC(\Sl)$ to $\cCO$
 step by step, each
is for   chord diagrams up to some degree. 

So we consider, for each fixed positive integer $n$, a sequence 
$T^n_m\in\cA(m), m=1,2,\dots$ with 
properties: (p1), (p2) (in which $\eta_m$ is replaced by $T^n_m$),
 and 
\begin{equation}  T^n_m=0 \quad\quad  \text{if $m<2n$}.\tag{p3}
\end{equation}

Certainly if the sequence $T^n_m, m=1,2, \dots$ satisfies these 
conditions, then the 
sequence $aT^n_m$ also
satisfies these conditions, for any  $a\in\cCO$.

It can be proved that up to this kind of ambiguity, there exists 
a unique sequence
$T^n_m\in\cA(m)/(R_{n+1})$, satisfying the above conditions. For 
the construction
of $T^n_m$, see \cite{LMO}. Actually, we will need only the 
element $T^n_{2n}$,
which had been described in \S\ref{iota}. The element 
$T^n_m\in\cA(m)$ is a linear combination
of uni-tri-valent graphs with $m$ external vertices and 
$m-2n$ internal vertices. Note that $T^n_{2n}$ is invariant under
any permutation of the $2n$ external vertices.

Now the mapping $\iota_n$ is defined as follows. 
For $\xi\in\cC(\Sl)$, a {\it chord diagram} with $m_i$ external 
vertices on the $i$-th solid component,  let us replace the 
$i$-th solid loop by
 $T^n_{m_i}$, then we get a linear combination 3-valent 
graphs which may contains
some loop components. Using relation $O_n$ (which says 
that a loop component
is equal to $-2n$), and deleting the part of degree $>n$, we obtain
$\iota_n(\xi)$ which is in $\Gr_{\le n}\cCO$. From this definition 
we get
\begin{pro} Suppose $\xi$ is a chord diagram with $k$ internal 
vertices, 
then $\iota_n(\xi)$ has at least $k$  vertices, hence has 
degree $\ge k/2$.
\end{pro}

\subsection{Finite type invariants for homology 3-spheres}

 We briefly
recall here some fundamental concepts, referring the reader to 
\cite{Ohtsuki,GO,GL}.

Let $\cM$ be the vector space over $\Bbb Q$ spanned by the set
 of all integral 
homology oriented 3-spheres.
Every invariant $I$ of homology 3-spheres with values in a {\it 
vector space} can
 be uniquely extended to a linear mapping on $\cM$.

For a framed link $L$ in $M$, let $M_L$ denote the 3-manifold 
obtained
 by performing
Dehn surgery on $L$. If $K$ is a formal linear combination 
$\sum_i c_iL_i$,
 where $L_i$
are framed links and $c_i$ are numbers, then $M_K$, by definition, is 
$\sum_ic_iM_{L_i}\in \cM$.

When $M$ is an oriented homology 3-sphere, there 
is a natural way
to identify the set of framings of a knot in $M$ with the set
 $\Bbb Z$ of integers; and we will use this identification.

A framed link $L$ in  $M$ is {\it unit-framed} if the framing of each 
component is $\pm1$;
$L$ is {\it algebraically split} if the linking number of every two 
components is 0.

Let $|L|$ be the number of components of $L$. Define
\begin{equation}\delta(L)=\sum_{L'\subset L}(-1)^{|L'|}L',
\label{del}
\end{equation}
where the sum is over all sublinks $L'\subset L$, including the
 empty link. 
 
Consider the following decreasing filtration in $\cM$. Let 
$\cF_n(\cM)$ be the vector space
 generated by
$M_{\delta(L)}$, where $M$ is an arbitrary homology 3-sphere and $L$ 
is a unit-framed
and algebraically split $n$-component link in $M$.
\begin{defn}(\cite{Ohtsuki}){\it 
An invariant $I$ of integral homology 3-spheres with values in a
 vector space is of order $\le n$ if $I(\cF_{n+1}(\cM))=0.$}
\end{defn}

 Ohtsuki showed  that $\cM/\cF_n(\cM)$ has finite 
dimension, 
see \cite{Ohtsuki}; this means the space of invariants of order 
less than $n$ is finite-dimensional.

An important result (see \cite{GL,GO}) is that $\cF_{3n}(\cM)=
\cF_{3n+1}(\cM)=\cF_{3n+2}(\cM)$.
Hence the set of invariants of order $3n$ is the same as the set of 
invariants of order $3n+2$.

Every invariant $I$ of order $3n$ with values in $\Bbb Q$ defines a 
linear 
form on $\cM/\cF_{3n+1}(\cM)$, hence restricts to a linear form on $\cF_{3n}(\cM)/\cF_{3n+1}(\cM)$. As in knot theory, it's important to
 understand the structure of $\cF_{3n}(\cM)/\cF_{3n+1}(\cM)$.

In \cite{GO}, a {\it surjective} linear mapping $\On:\Gr_n(\cCO)
\to\cF_{3n}/\cF_{3n+1}(\cM)$ was constructed. For a 
definition of $\On$, see \S\ref{sec:lemma1} below.
It follows that if $I$ is an invariant of order $3n$, then by 
combining
with $\On$, $I$ defines a linear mapping
from $\Gr_n(\cCO)$ to $\Bbb Q$, called the {\it corresponding linear 
form} of $I$. Note that if the corresponding linear form is 0,
then the invariant is of order $< 3n$.

Is every linear
mapping from $\Gr_n(\cCO)$ to $\Bbb Q$ the linear form of some
invariant of order $\le 3n$? (Question 1 in \cite{GO}).
We will see that the answer is positive.
This is equivalent to the fact that the above mapping $\On$ 
is an isomorphism
 of vector spaces.

\section{Main results}
\subsection{Universality of $\Omega$}
We will show the following 

\begin{thm}{\rm (Main Theorem)}

a) For homology 3-spheres,
$\Omega_n$ is an invariant of order $\le 3n$.

b) The mapping $\On$ is an isomorphism, i.e. every linear form 
from $\Gr_n(\cCO)$
 to $\Bbb Q$ is a linear form of some finite invariant of order $3n$.

c) $\Omega_n$ is a universal invariant of order $3n$, i.e. if \/ 
$\Omega_n(M)=\Omega_n(M')$
then $I(M)=I(M')$ for every invariant $I$ of order $\le 3n$.
\end{thm} 
We need the following two lemmas, 
whose 
proofs will be presented in Section \ref{proofs}.
\begin{lem} Suppose that $\Gamma$ is a 3-valent vertex-oriented 
graph in $\cCO$ 
of degree $n$. Then for any representative $\alpha(\Gamma)$   of
 $\On(\Gamma)$ in $\cF_{3n}(\cM)$, one has
\begin{equation}\Omega_n(\alpha(\Gamma))=(-1)^n\Gamma +\Omega_n(M_1),
\label{equ:lemma1}\end{equation}
where $M_1$ is in $\cF_{3n+3}(\cM)$.
\label{lemma1}\end{lem}
\begin{lem}\label{lemma2}
If $N\in\cF_{6n+1}(\cM)$, then $\Omega_n(N)=0$.\end{lem}
This lemma is certainly weaker than part $a)$ of Main Theorem.
\begin{pf}[of Main Theorem]

a) Suppose $N\in\cF_{3n+1}(\cM)=\cF_{3n+3}(\cM)$ we have to show that
$\Omega_n(N)=0$.

Since $O^*_{3n+3}:\Gr_{n+1}(\cCO)\to\cF_{3n+3}(\cM)/\cF_{3n+6}(\cM)$ 
is surjective, 
there is $\Gamma_1\in\Gr_{n+1}(\cCO)$ such that 
$$N-\alpha(\Gamma_1)=N_1\in \cF_{3n+6}(\cM),$$
 where $\alpha(\Gamma_1)$ 
is
any representative of $O^*_{3n+3}(\Gamma_1)$ in $\cF_{3n+3}(\cM)$.
Apply $\Omega_{n+1}$ to both sides, using Lemma \ref{lemma1} with 
$n$
 replaced by $n+1$; we see that, for some  $M_1\in\cF_{3n+6}(\cM)$,
$$\Omega_{n+1}(N)=\pm\Gamma_1 + \Omega_{n+1}(N_1) + \Omega_{n+1}(M_1)
=\pm\Gamma_1+\Omega(N_1+M_1).$$
 Since $\Omega_n=\Gr_{\le n}\Omega_{n+1}$ 
(see Proposition
\ref{ax}) and $\Gamma_1$ is of degree $n+1$, one has
$$\Omega_n(N)=\Omega_n(M_1+N_1).$$
Here $(M_1+N_1)$ is in $\cF_{3n+6}$, while $N$ is in 
$\cF_{3n+3}(\cM)$. Repeat this argument
 $n$ times, one finds $N'$ in $\cF_{6m+3}(\cM)$ such that
$$\Omega_n(N)=\Omega_n{N'}.$$
  Now Lemma \ref{lemma2} says 
that 
$\Omega_n(N')=0$.
This completes the proof of $a)$.

b) We will show an inverse of $\On:\Gr_{n}(\cCO)\to\cF_{3n}
(\cM)/\cF_{3n+3}(\cM)$.
 Let us consider equation (\ref{equ:lemma1}) again. Since $M_1$ 
in the equation
is in $\cF_{3n+3}(\cM)$, by part $a)$, we have 
$\Omega_n(M_1)=0$. Hence
$$\Omega_n(\alpha(\Gamma))=(-1)^n\Gamma.$$
This is true for every representative $\alpha(\Gamma)$ of 
$\On(\Gamma)$ 
in $\cF_{3n}(\cM)$, hence $(-1)^n\Omega_n$ restricts to a
 well-defined linear
 mapping from $\cF_{3n}(\cM)/\cF_{3n+3}(\cM)$
to $\Gr_n{\cCO}$, which, by the above identity, is inverse to $\On$.

c) We use induction on $n$.  Suppose $\Omega_n(M)=\Omega_n(M')\in 
\Gr_{\le n}(\cCO)$ and $I$ is an invariant 
of order $n$ with values in $\Bbb Q$. Let $W:\Gr_{n}(\cCO)\to 
\Bbb Q$ be the corresponding
linear form of $I$. Let $I_1=W(\Gr_n(\Omega_n))$. Then $I_1$ is an 
invariant 
of order $\le 3n$ and $I_1(M)=I_1(M')$. 

Note that $I-I_1$ is an 
invariant 
of order $\le 3n$. It is not hard to verify that 
the corresponding linear form of $I-I_1$ from  $\Gr_{n}(\cCO)$
 to $\Bbb Q$ is 0,
 hence $I-I_1$ is of
order $\le (3n-1)$. By induction, $(I-I_1)(M)=(I-I_1)(M')$. Hence 
$I(M)=I(M')$.
\end{pf}

Let the {\it product} of two homology 3-spheres be their connected 
sum. The unit is $S^3$.
Define a co-product: $\Delta (M)=M\otimes M$,  and 
co-unit: $\varepsilon(M)=1$, for every homology 3-sphere $M$. 
Then $\cM$ becomes
a commutative co-commutative Hopf algebra.

It is easy to see that the product is compatible with the filtration
$$\cM=\cF_0{\cM}\supset\cF_3{\cM}\supset\dots\supset\cF_{3n}(\cM)
\dots$$

This means, if $N_1\in\cF_{n_1}(\cM)$ and $N_2\in\cF_{n_2}(\cM)$, 
then their product is in $\cF_{n_1+n_2}(\cM)$.

 Hence the algebra structure of $\cM$ induces an algebra structure
 on the
 associated (complete) graded vector space 
$${\cal {GM}}=\prod _{i=0}^{\infty}\cF_{3i}(\cM)/\cF_{3i+3}(\cM).$$

A nice interpretation of the main theorem is the following.
\begin{thm}
\label{main2}
a) The mapping $\Omega: \cM\to \cCO$ is a {\em Hopf algebra} 
homomorphism
which maps $\cF_{3n}(\cM)$ to $\Gr_{\ge n}\cCO$.

b) $\Omega$ induces an {\em algebra} isomorphism between $\cal {GM}$
and $\cCO$.
\end{thm}

This theorem follows immediately from the previous one and 
Propositions \ref{grouplike}, \ref{connectedsum}. Part $b)$ shows
that $\cal {GM}$ has a structure of a commutative co-commutative
Hopf algebra.

\subsection{Some further properties of $\Omega$}

\begin{pro}\label{x40}
a)
$\Omega_n$ is of order exactly $3n$.

b)  For every $n$ there is {\em an integral homology
3-sphere} $M_n$ such that $\Omega(M_n)=1+\text{degree} \ge n$, and the 
$n$-th degree part is not 0.
\end{pro}
\begin{pf} 
a) This is equivalent to the fact that $\Gr_n(\cCO)$ has positive 
dimension. It's enough
to show that there is a non-zero linear form on $\Gr_n(\cCO)$. One 
can use the weight
system coming from simple Lie algebra, as defined in 
\cite{Ko1,BN1,LM2},
 to define
non-zero linear form. The details are left for the reader. $a)$ also
 follows from $b)$.

b)  Consider a solid loop and a dashed loop, each has $2n$ marked 
points $x_1,x_2,\dots,x_{2n}$ and $y_1,y_2,\dots,y_{2n}$, counting
 counterclockwise.
 Connect $x_i$
with $y_i$ by a dashed line, for every $i=1,2,\dots,2n$, we get a 
chord 
diagram $\xi_n$ on a solid loop. 
In \cite{Ng}, a knot $K$ was constructed with property that for 
every 
knot invariant  $I$ of degree $<2n$, we have 
$I(K)=I(\text{unknot})$ and for 
every invariant  $I$ of degree
$2n$ we have $I(K)=I(\xi_n)$. It follows that 
if  $K$ has 
framing 0, then 
$$\hZ(K)=\nu(1+\xi_n) + \text{(elements of degree $>2n$)}.$$
Hence if $K'$ is the knot $K$ with framing 1, we have
$$\cZ(K')=\nu^2e^{\theta/2}(1+\xi_n)+\text{(elements of 
degree $>2n$)}.$$
Here $\theta\in\cC(S^1)$ is the chord diagram of degree 1 
with one isolated
chord; $e^{\theta/2}$ is the contribution of framing 1 
(see \cite{LM2}, Theorem 3).
Let $M_n$ be obtained by Dehn surgery on $K'$, then, by definition
$$\Omega_n(M_n)=\frac{\iota_n(\cZ(K'))}{\iota_n(U_+)}=\frac
{\iota_n[\nu^2e^{\theta/2}(1+\xi_n)]}{\iota_n(\nu^2e^{\theta/2})}.$$
Since $\xi_n$ is of degree $2n$, $\iota_n$ annihilates any 
chord diagram
in $\cC(S^1)$ 
 of degree $> 2n$, and $\iota_n(U_+)=\iota_n(\nu^2e^{\theta/2})=
(-1)^n+...$, we see that 
$\Omega_n(M_n)=1+(-1)^n\iota_n(\xi_n)$.

Using the definition of $\iota_n$ (the one with $T_{2n}^n$ 
described in 
\S\ref{uni-inva}), combining with $sl_2$-weight system, one can 
 prove that $\iota_n(\xi_n)\not=0$.
\end{pf}

In \cite{Lin}, X. S.  Lin asked whether there exists an operation,
 similar to 
the Stanford one in the knot case, on homology 3-spheres which 
does not alter the values of 
any invariant of order $\le n$. We now present such an operation.

Suppose $M$ is obtained by Dehn surgery on a unit-framed 
 algebraically split link $L$.
Suppose that $T$ is a part of $L$ which is a trivial tangle 
with $m$ strands. 
Let $G_j(P_m), j=0,1,2,\dots$ be the lower central series of 
the pure braid
group $P_m$ on $m$ strands, defined by $G_0(P_m)=P_m,G_{j+1}
(P_m)=[G_j(P_m),P_m]$.
Replacing $T$
by an element $\gamma$ in the $(2n+2)$-th group $G_{2n+2}(P_m)$, 
we obtain a unit-framed algebraically split link $L'$. 
The homology 3-sphere $S^3_{L'}$ is said
{\it to be obtained from $M$ by an $s_n$ operation}. 
\begin{pro}
If $M'$ is obtained from $M$ by an $s_n$ operation, then 
$\Omega_n(M)=\Omega_n(M')$. 
Hence the values of any invariant of order $\le 3n$ on $M$ 
and on $M'$ are the same.
\end{pro}
\begin{pf}
We first consider $\hZ(\gamma)\in\cP_m$.
By a result of Stanford \cite{Sta}
(see also \cite{Kohno}), the values of any pure braid 
invariant of order $\le (2n+1)$ 
on $T$ (the trivial tangle)  and $\gamma$ are the same. Hence 
$$\Gr_{\le (2n+1)}\hZ(\gamma)=1.$$
From Theorem 4.2 of \cite{LM3}, it follows that $\hZ(\gamma)$
is a group-like element, i.e. $\hZ(\gamma)=\exp(\xi)$,
where $\xi$ is primitive. The above identity shows that $\xi$ 
is a linear combination of chord diagrams of degree $\ge (2n+2)$.
Hence from Proposition \ref{ifilter} one sees that 
$$\hZ(\gamma)=
1+ (\text{elements of i-filter $\ge (2n+1)$}).$$ 
Since $\iota_n$ annihilates any element of i-filter 
$\ge (2n+1)$, we have that $\Omega_n(M)=\Omega_n(M')$.
\end{pf}
Note that we proved this fact only for invariants with {\it 
rational } values. 
It would be interesting to give a direct geometric proof of this 
fact, and to establish the same fact
for invariants with values in any {\it abelian} group.

From this proposition, it is easy to construct, for any given homology
3-sphere $M$ and any positive integer $n$, infinitely many 
homology 3-spheres which are indistinguishable from $M$ by all
 invariants of order $\le n$.

\subsection{Some other corollaries; remarks}
X. S.  Lin proved (see \cite{Lin})
\begin{pro}
The space of finite type invariants of homology 3-spheres with
 rational values is a commutative co-commutative Hopf algebra.
\end{pro}
 Here this fact is a consequence of Theorem \ref{main2}, $b)$,
as observed in \cite{GO}, Remark 1.10.

 Suppose $I$ is an invariant of integral homology 3-spheres,
$K$ is a knot in $S^3$. Let us define $\lambda_I(K)=I(S^3_K)$, 
where we provide $K$ with framing 1. 

\begin{pro} If $I$ is a homology 3-sphere invariant of order $3n$, 
then $\lambda_I$ is a knot invariant of order $2n$.
\end{pro}

This had been 
conjectured by Garoufalidis, and proved by Habegger, see \cite{Hab}. 
Here this fact follows easily from the main
theorem: for a  knot $K$, the universal invariant $\Omega_n(S^3_K)$ 
is computed using $\Gr_{\le 2n}(\hZ(K))$, hence $I_{\Omega_n}(K)$ is
 derived
from $\Gr_{\le 2n}\hZ(K)$ which is a knot invariant of order $\le 2n$.

\begin{rem}

1. The main theorem can be reformulated and proved (in a similar way) 
for
 {\it rational}
homology 3-spheres, see \cite{GO2} for theory of finite type 
invariants 
of rational homology 3-spheres. The details will appear elsewhere.

2. The construction of $\Omega$ suggests that the relation between 
$\Omega$, 
combined with $sl_2$ weight, and $sl_2$ quantum invariant at 
roots of unity
maybe expressed through Ohtsuki polynomial in \cite{Oht2}, or 
the invariant
defined in \cite{Oht3}.
\end{rem}

\subsection{ Groups of homology 3-spheres}
In analogy with the knot case (see \cite{Gus,Ng,SN}), we say that two
integral homology spheres $M_1,M_2$ are $V_n$-{\it equivalent} if 
they are indistinguishable by any {\it rational} invariant of order 
$\le n$; or the same as: $(M_1-M_2)\in\cF_{3n+1}(\cM)$, or
$\Omega_n(M_1)=\Omega_n(M_2)$.

Let $\VM_n$ be the set of all homology 3-spheres, regarded up to
the $V_n$-equivalence relation. This set $\VM_n$ is a semigroup, 
where the product is the connected sum.

\begin{thm} \label{existence}

For every {\em connected} 3-valent vertex-oriented graph $\Gamma$
of degree $n$,
there are homology 3-spheres $M^{\pm}(\Gamma)$ such that
$$\Omega(M^{\pm}(\Gamma))=1+\pm\Gamma + (\text{\rm elements
 of degree $>n$}).$$
Moreover, $M^{\pm}(\Gamma)$ are obtained by Dehn surgery on links
of $n$ components.
\end{thm}
This is stronger than proposition \ref{x40}.
From Theorem \ref{existence}, one gets

\begin{thm}\label{group}

The semi-group $\VM_n$ is a {\em group}. This means, for every
homology 3-sphere $M$, there is another homology 3-spheres $M'$
such that $M\#M'$ is $V_n$-equivalent to ~0.

The group $\VM_n$ is free abelian of rank equal to the dimension
of the subspace of $\Grn(\cCO)$, spanned by {\em connected}
3-valent vertex-oriented graphs.
\end{thm}
This answers the first part of question 2 in \cite{Lin}.

Proofs of Theorems \ref{existence} and \ref{group}  will 
appear elsewhere.

\section{Proofs}\label{proofs}
\subsection{The mapping $O^*_{3n}$}
\label{sec:lemma1}

Let's consider a 3-valent vertex-oriented graph $\Gamma$ of
 degree $n$ in
 $\bR^3$. Like with links, we can represent $\Gamma$ using
a generic projection on the plane with decoration at every 
double point to
 indicate
over/under-crossing. We suppose that $\Gamma$ is equal to this
 projection 
everywhere
except for a neighborhood of the double points, that 
 the orientation at every vertex is given by counterclockwise 
direction of
 the plane, and
that the three edges incident to a 3-valent vertex is going 
downwards in a 
neighborhood of the vertex as in Figure \ref{3-valent}. Here 
edges of $\Gamma$ are depicted by solid lines.

$\Gamma$ has $2n$ trivalent vertices and $3n$ edges. Using 
$2n$ pairs  of 
horizontal lines on the plane,
each pair consists of a line above and a line below a vertex, 
we can 
decompose $\Gamma$ as:
\begin{equation}\Gamma=T_1V_1T_2V_2\dots T_{2n}V_{2n}T_{2n+1}.
\label{decomposition}
\end{equation}
Here $T_i$ are some non-oriented tangles while $V_i$ are trivial 
tangles with a 3-valent vertex part as in Figure \ref{3-valent};
and product in the right hand side  means ``placing on top".
\begin{figure}[htp]
\centerline{
\begin{picture}(120,60)
\thicklines \put(0,60){\line(0,-1){60}}
\put(20,60){\line(0,-1){60}}
\put(120,60){\line(0,-1){60}}
\put(70,30){\line(0,-1){30}}
\put(70,30){\line(1,-1){30}}
\put(70,30){\line(-1,-1){30}}
\end{picture}
}
\caption{$V_i$\label{3-valent}}
\end{figure}

Now we define $\beta(T_i)$ and $\beta(V_i)$ as follows. Let 
$\beta(T_i)$ be 
the {\it tangle} obtained by replacing each component of $T_i$ 
by a pair
of parallel push-offs (on the plane) of this component. 
The orientation on
$\beta(T_i)$ is chosen so that at  boundary points of 
a pair of push-offs, the left boundary point is pointing 
upwards, and the right boundary point  downwards.

\begin{figure}[htp]
\centerline{
\begin{picture}(120,60)
\thicklines \put(0,60){\line(0,-1){60}}
\put(3,60){\line(0,-1){60}}
\put(20,60){\line(0,-1){60}}
\put(23,60){\line(0,-1){60}}
\put(120,60){\line(0,-1){60}}
\put(123,60){\line(0,-1){60}}
\put(36,37){\line(0,-1){37}}
\put(55,37){\oval(38,14)[t]}
\put(73,33){\line(-1,-1){16}}
\put(53,13){\line(-1,-1){13}}
\put(70,0){\line(-1,1){17}}
\put(63,17){\oval(20,37)[tl]}
\put(69,14){\oval(38,43)[tr]}
\put(84,10){\line(-1,-1){10}}
\put(100,0){\line(-6,5){31}}
\put(69,30){\oval(8,24)[tl]}
\put(75,0){\oval(56,84)[tr]}
\end{picture}
\hskip 1cm
\begin{picture}(10,60)
\put(5,30){\makebox(0,0){$-$}}
\end{picture}
\hskip 1cm
\begin{picture}(120,60)
\thicklines \put(0,60){\line(0,-1){60}}
\put(3,60){\line(0,-1){60}}
\put(20,60){\line(0,-1){60}}
\put(23,60){\line(0,-1){60}}
\put(120,60){\line(0,-1){60}}
\put(123,60){\line(0,-1){60}}
\put(44,0){\oval(4,40)[t]}
\put(72,0){\oval(4,40)[t]}
\put(100,0){\oval(4,40)[t]}
\end{picture}
}
\caption{$\beta(V_i)$\label{3valent}}
\end{figure}
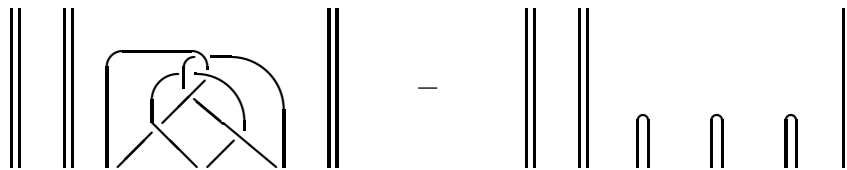

 $\beta(V_i)$ is not a {\it single tangle}, but rather  the
{\it formal difference of two tangles}, $u-v$. The first tangle 
$u$ is obtained from $V_i$
by replacing each vertical
component of $V_i$ by a pair of its parallel push-offs and 
replacing the 
3-valent
vertex by 3 arcs linked with each other as in the Borromean 
ring; in the 
second tangle $v$, the 3-valent vertex is replaced by 3 arcs 
which do not link with 
each other. For example, the
$V_i$ of Figure \ref{3-valent} has $\beta(V_i)$ as in 
Figure~\ref{3valent}. The orientation is chosen so that at boundary
points of a pairs of push-offs, or at boundary points of a non-vertical
arc, the left point is pointing upwards, and the right point is pointing
downwards.

$\beta(T_i)$ will be regarded as a {\it framed} oriented tangle, while 
$\beta(V_i)$  will be
regarded as the difference of two {\it framed} oriented tangles, where 
the framing vector
at every point is given by the  vector
perpendicular to the page and pointing at us.

We define 
$$\beta(\Gamma)=\beta(T_1)\beta(V_1)\dots\beta(V_{2n})
\beta(T_{2n+1}).$$
This $\beta(\Gamma)$ is not a single link, but a {\it 
formal linear combination}
of $2^{2n}$ algebraically split framed oriented links, each has 
$3n$ components
 with framing 0. 

Changing the framing of each link component to 1, from 
$\beta(\Gamma)$ we get
$\tilde{\beta}(\Gamma)$ which is a linear combination of
 $2^{2n}$ unit-framed
algebraically split oriented links. 

Then $\alpha(\Gamma)=
S^3_{\delta(\tilde{\beta}(\Gamma))}$
 is an element in 
$\cF_{3n}(\cM)$. Recall that the operator $\delta$ is given by 
(\ref{del}). Define
$\On(\Gamma)$ as the image of $\alpha(\Gamma)$ in 
$\cF_{3n}(\cM)/\cF_{3n+1}(\cM)$.
In \cite{GO} it was proved that $\On$ is a well-defined 
{\it surjective} 
linear mapping from $\Grn(\cCO)$ to $\cF_{3n}
(\cM)/\cF_{3n+1}(\cM)$.
\begin{lem} One has  $\delta(\tilde{\beta}(\Gamma))
=\tilde{\beta}(\Gamma)$.
Hence $\alpha(\Gamma)=S^3_{\tilde{\beta}(\Gamma)}$.
\label{phu11} \end{lem}
This  combinatorial lemma follows easily from the definition and
 a simple induction.

\subsection{Proof of Lemma \protect{\ref{lemma1}}} We will prove
\begin{lem} For every 3-valent vertex-oriented graph $\Gamma$ 
in $\Grn(\cCO)$ one has:
\begin{equation}\Omega_n(S^3_{\tilde{\beta}(\Gamma)})=(-1)^n\Gamma.
\label{main}
\end{equation}
\label{xphu}
\end{lem}
Note that Lemma \ref{lemma1} follows from this lemma, since $\alpha(\Gamma)=S^3_{\tilde{\beta}(\Gamma)}$ (by Lemma 
\ref{phu11}) and any two  representatives
of $\On(\Gamma)$ differ by an element in $\cF_{3n+3}(\cM)$.

For a formal linear combination $\sum c_iL_i$ of framed oriented links,
 we define 
$\cZ(\sum c_iL_i)$ as $\sum c_i\cZ(L_i)$. 
Then $\cZ(\tilde{\beta}(\Gamma))
\in\cC(\sqcup^{3n}S^1)$.
The left hand side of (\ref{main}), by definition, is
$$\Omega_n(S^3_{\tilde{\beta}(\Gamma)})=
\frac{\iota_n[\cZ(\tilde{\beta}(\Gamma))]}
{\iota_n(U_+)^{-3n}}.$$
Since $\iota_n(U_+)$ has the form $(-1)^n+$ 
elements of degree $\ge 1$ 
(see \cite{LMO}), in order to prove (\ref{main}) 
it suffices to show that
\begin{equation}\iota_n[\cZ(\tilde{\beta}(\Gamma))]
=\Gamma.\label{main1} 
\end{equation}

For $V_i$ in the decomposition (\ref{decomposition}) of $\Gamma$, 
$\beta(V_i)=u-v$, where
$u-v$ are tangles in Figure \ref{3valent}. Note that as abstract 
1-manifolds, $u$ and $v$ are homeomorphic, and $\hZ(\beta(V_i))=
\hZ(u)-\hZ(v)$ is in $\cC(v)$. Let $\xi(V_i)$ be the {\it chord
 diagram} in $\cC(v)$ with exactly one dashed component with 
1 internal
vertex  which is 
connected to the 3 non-vertical solid lines of $v$, as 
Figure~\ref{xi(Vi)}.
Then $\xi(V_i)$ has {\rm i}-filter 1.
\begin{figure}[htp]
\centerline{
\begin{picture}(120,60)
\put(44,40){\dashbox{2}(56,0){}}
\put(44,20){\dashbox{2}(0,20){}}
\put(100,20){\dashbox{2}(0,20){}}
\put(72,20){\dashbox{2}(0,20){}}
\thicklines \put(0,60){\line(0,-1){60}}
\put(3,60){\line(0,-1){60}}
\put(20,60){\line(0,-1){60}}
\put(23,60){\line(0,-1){60}}
\put(120,60){\line(0,-1){60}}
\put(123,60){\line(0,-1){60}}
\put(44,0){\oval(4,40)[t]}
\put(72,0){\oval(4,40)[t]}
\put(100,0){\oval(4,40)[t]}
\put(72,40){\circle*{2.5}}
\end{picture}
}
\caption{\label{xi(Vi)} $\xi(V_i)$}
\end{figure}

 \begin{lem} One has
 $\hZ(\beta(V_i))=\hZ(u)-\hZ(v)=\xi(V_i)+ 
(\text{elements of {\rm i}-filter $\ge 2$})$.
\label{phu111}
\end{lem}

A proof of the lemma is given in \S\ref{x50}, after a proof
of a similar fact.

Let $\xi(\Gamma)$ be the {\it chord diagram} with support $3n$ solid 
loops obtained
as follows: 
$$\xi(\Gamma)=T_1\,\,\xi(V_1)\,\,T_2\,\,\xi(V_2)\dots T_{2n}
\,\,\xi(V_{2n})\,\,T_{2n+1}.$$

On every solid loop of $\xi(\Gamma)$ there are exactly two
 external vertices.
It is important to notice that  if we remove every solid loop of 
$\xi(\Gamma)$ and connect the two external vertices of each solid loop 
by a dashed line, then we get exactly $\Gamma$.

Now we consider $$\hZ(\beta(\Gamma))=\hZ(\beta(T_1))\,
\hZ(\beta(V_1))\dots
\hZ(\beta(V_{2n}))\,
\hZ(\beta(T_{2n+1})).$$
 Recall that  $\Omega_n$ annihilates any chord diagram of 
{\rm i}-filter $>2n$.

Using Lemma \ref{phu111} we see that $\hZ(\beta(\Gamma))$ 
is, modulo 
elements of {\rm i}-filter $>2n$, a linear combination 
of chord diagrams, 
each contains $\xi(\Gamma)$ as a subdiagram,
i.e. each is $\xi(\Gamma)$ plus some extra chords. The extra
 chords are 
contributions from $\hZ(\beta(T_i))$. Note that every 
$\beta(T_i)$ is obtained
 by first taking
parallel push-offs of every component of  the tangle 
$T_i$, then putting 
orientations on components such that
the two push-offs are of opposite orientation. Hence 
Theorem \ref{theorem1}
  says that 
$\hZ(\beta(T_i))$ is $1+$ a linear combination  of terms 
of the form 
$\xi-\xi'$, where $\xi$ and $\xi'$ are chord diagrams 
(in $\cC[\beta(T_i)]$) 
identical everywhere, except for a ball in which they are as in Figure 
\ref{difference}.

\begin{figure}[htp]
\centerline{
\begin{picture}(120,60)
\put(-15,40){\dashbox{2}(15,0){}}
\put(100,40){\dashbox{2}(20,0){}}
\put(0,5){\makebox(0,0){$\xi$}}
\put(115,5){\makebox(0,0){$\xi'$}}
\thicklines \put(0,20){\vector(0,1){40}}
\put(5,60){\vector(0,-1){40}}
\put(120,60){\vector(0,-1){40}}
\put(115,20){\vector(0,1){40}}
\end{picture}
}
\caption{\label{difference}}\end{figure}

The two solid lines here are of the same component of the 
support of chord diagrams
in $\cC(\sqcup^{3n}S^1)$. Using the  STU relation
we can always move any external vertex around the solid loop; when  
exchanging two external vertices, we have to
add a chord diagram with an extra {\it internal} vertex. Hence we
 can cancel $\xi$ and $\xi'$, and get
\begin{equation}
\hZ(\beta(\Gamma))=\xi(\Gamma)+(\text{elements of 
{\rm i}-filter $>2n$}).
\end{equation}

Now we prove (\ref{main1}).
The left hand side of (\ref{main1}) is, by the definitions
 of $\cZ$ and $\tilde{\beta}$, 
$$\iota_n[\cZ(\tilde{\beta}(\Gamma))]=
\iota_n[(\nu e^{\theta/2})^{\otimes 3n}\# \hZ(\beta(\Gamma))],$$
where $\theta$ is the chord diagram on a solid loop with exactly one 
isolated chord; $e^{\theta/2}$ is the contribution of the framing 1 
(see Theorem 3 of \cite{LM2}). 

Note that $\xi(\Gamma)$ already has $2n$ internal vertices and 
$\nu$ has the form $1+$ ({\rm i}-filter $\ge 1$), hence
we can delete
all the $\nu$'s in the formula.

  Every  solid component of $\xi(\Gamma)$
 has 2 external vertices.
Since  $\iota_n$ 
annihilates any chord diagram
with less than $2n$ external vertices on a solid component, 
we see that 
\begin{equation}\label{phu}
\iota_n[(\nu e^{\theta/2})^{\otimes 3n} \#\hZ(\beta(\Gamma))]= 
\iota_n\{[\frac{1}
{2^{n-1}(n-1)!}]^{3n}\tilde {\xi}(\Gamma)\},\end{equation}
where $\tilde{\xi}(\Gamma)$ is obtained from $\xi(\Gamma)$ by
 adding $(n-1)$ 
isolated chords
to each solid component. The isolated chords and the
 coefficients 
$\frac{1}{2^{n-1}(n-1)!}$
are contributions of $e^{\theta/2}$.

$\tilde{\xi}(\Gamma)$ is a chord diagram with $2n$ external 
vertices 
on each solid loops,
$(2n-2)$ of them are external vertices of the $(n-1)$ 
isolated chords. 

\begin{lem} One has
$$\iota_n(
\text{
\begin{picture}(80,15)(0,5.5)
\put(27,0){\dashbox{2}(0,15){}}
\put(30,0){\dashbox{2}(0,15){}}
\put(51,0){\dashbox{2}(0,15){}}
\put(41,7.5){\makebox(0,0){$\cdots$}}
\put(0,7.5){\dashbox{2}(15,0){}}
\put(65,7.5){\dashbox{2}(15,0){}}
\thicklines
\put(40,7.5){\oval(50,15)}\end{picture}}\, )=(-2)^{n-1}(n-1)!\,\,(\text{
\begin{picture}(34,5)
\put(2,3){\dashbox{2}(30,0){}}
\end{picture}} \,).$$
Here the chord diagram in the left hand side has $(n-1)$ 
vertical chords.
The result holds true if these vertical dashed chords are 
replaced by 
isolated chords.
\label{x60}
\end{lem}
\begin{pf} 
By definition, we have to replace the solid loop by element 
$T^n_{2n}$
 (see \S\ref{uni-inva}), with
the convention that a dashed loop is equal $-2n$, by relation $O_n$.
 Because
of the symmetry of $T^n_{2n}$, the left hand side does not depend 
on whether
the $(n-1)$ dashed lines are vertical chords or isolated chords.
The proof of the Lemma now follows by an induction.
\end{pf}

Note that if we remove every solid loop of $\xi(\Gamma)$ and connect
the two external vertices of each solid loop be a dashed line, then we
get exactly $\Gamma$.
Hence from Lemma \ref{x60} 
we get
$$\iota_n(\tilde {\xi}(\Gamma))=[(-2)^{n-1}(n-1)!]^{3n}\,\, \Gamma.$$

This, together with (\ref{phu}) proves (\ref{main1}), and hence
Lemma \ref{lemma1}.

\subsection{Beginning of Proof of Lemma \protect{\ref{lemma2}}}

 Let $K$ be a unit-framed algebraically split link in a homology 
3-sphere $M$ with $(6n+1)$ 
components. We have to show that $\Omega_n(M_{\delta(K)})=0$. 
The following 
is well-known and can be proved using the Kirby calculus for links in 
3-manifolds and
Lemma 1.1 of \cite{Hos}.

\begin{lem} There is a unit-framed and algebraically split link 
 $L$ in $S^3$ which is a union 
of two links $L_1,L_2$ 
such that after Dehn surgery on $L_1$, from $S^3$ we get $M$ 
and from $L_2$
we get $K$. ($L_1$ and $L_2$ may link with each other in $S^3$).
\end{lem}

Let 
$$\tilde \delta_{L_2}(L)=\sum_{L_2'\subset 
L_2}(-1)^{|L_2'|}(L_1,L_2''),$$ 
where the sum
 is over the set of all sublinks
$L_2'$ of $L_2$; and $(L_1,L_2'')$ is 
the link
 obtained from $L=(L_1,L_2)$ by replacing every component of 
$L_2\setminus 
L_2'$ by a trivial knot with the same framing; the trivial knots are 
not linked with each other nor with $L$.

Let  $(L_1,L_2')$ be
obtained from $(L_1,L_2)$ by replacing $L_2$ with $L_2'$.

By definition,
$$\Omega_n(M_{\delta(L)})=\sum_{L_2'\subset L_2}(-1)^{|L_2'|}\Omega_n
[S^3_{(L_1,L_2')}]=
\sum_{L_2'\subset L_2}(-1)^{|L_2'|}
\frac{\iota_n(\cZ(L_1,L_2'))}{\iota_n(U_+)^{-\sigma_+(L_1,L_2')}
\iota_n(U_-)^{-\sigma_-(L_1,L_2')}}.$$

Note that $U_\pm$ are the trivial knots with framing $\pm1$, and that 
$\iota_n(a\sqcup b)=\iota_n(a)\iota_n(b)$ for two links $a,b$. Using 
the same 
denominator
$ \iota_n(U_+)^{-\sigma_+}\iota_n(U_-)^{-\sigma_-}$, where 
$\sigma_{\pm}=
\sigma_{\pm}(L_1,L_2)$, for all the terms of the last sum, we 
arrive at the following important formula:
$$\Omega_n(M_{\delta(K)})=\frac{\iota_n[\cZ 
(\tilde{\delta}_{L_2}(L_1,L_2))]}{\iota_n
(U_+)^{-\sigma_+}\iota_n(U_-)^{-\sigma_-}}.$$

Hence in order to prove Lemma \ref{lemma2}, it's sufficient to prove 
that
\begin{equation}
\iota_n[\cZ(\tilde {\delta}_{L_2}(L_1,L_2))]=0
\label{equ:lemma2}
\end{equation}
Here $(L_1,L_2)$ is a unit-framed algebraically split link in 
$S^3$ with $|L_2|=6n+1$.

The proof of (\ref{equ:lemma2}) may look complicated, but the idea 
is simple. The point is that $\tilde{\delta}_{L_2}(L)$ is a repeated
difference formula. The simplest case is when $|L_2|=1$ 
(only one step). In this case 
$\tilde{\delta}_{L_2}(L)=L-L'$, where $L'$ is obtained from $L$ by
replacing the only component of $L_2$ by a trivial knot with the same
framing. Then it can be proved that $\cZ(L-L')$ is of i-filter $\ge 1$.
(The reason for this phenomenon is that the linking number
of $L_2$ and other components is 0).

When $|L_2|>1$ we have a repeated difference sum for 
$\cZ(\tilde {\delta}_{L_2}(L_1,L_2))$; we expect that
$\cZ(\tilde {\delta}_{L_2}(L_1,L_2))$ has i-filter $\ge (2n+1)$,
if $|L_2|$ is sufficiently large ($|L_2|=6n+1$ is enough). 
It follows then $\iota_n[\cZ(\tilde {\delta}_{L_2}(L_1,L_2))]=0$.
A rigorous proof requires a lot of technical stuff.

\subsection{Some preparations}
Let $Y$ be a submanifold of $X$; $Y$ consists of several 
components of $X$. 
We say that an element $\xi\in\cC(X)$ is {\it {\rm i}-near} 
$Y$ if $\xi$ is a linear combination of chord diagrams, 
each has  internal vertices {\it near} every component of $Y$. 
Here an internal
vertex is near a solid component if the vertex is connected to 
the solid 
component by a dashed line of the dashed graph.

Recall that $\cP_m$, the space of chord diagrams on $m$ numbered 
solid lines,
is a co-commutative Hopf algebra. The space of primitive elements
 is spanned 
by chord diagrams with {\it connected} dashed graph.

A tangle  consisting of $m$ strings (no loops) such that the two
end points of each string project vertically to the same point 
is called a {\it string link}. Suppose 
that $T$ is an $m$-component framed
string link whose orientation at every boundary points is 
downward. A corollary of Theorem 4.2 of \cite{LM3} is that $\hZ(T)\in 
\cP_m$ is group-like.
Hence $\hZ(T)=\exp(\xi)$, where $\xi$ is primitive.

Suppose $T'$ is obtained from $T$ by removing a component $C$.
Let 
$\varepsilon_C:\cP_m\to\cP_{m-1}$ be the linear mapping which send 
a chord diagram $\xi\in\cP_m$
to $\varepsilon_C(\xi)$ obtained from $\xi$ by removing the 
 solid line $C$ 
if there is no external vertex on $C$; otherwise 
$\varepsilon_C(\xi)=0$.

The relation between $\hZ(T)$ and $\hZ(T')$ is expressed by 
(see \cite{LM1,LM2}):
\begin{equation}
\hZ(T')=\varepsilon_C[\hZ(T)].
\label{x15}
\end{equation}

\begin{figure}[htp]
\centerline{
\begin{picture}(120,110)(0,40)
\thicklines
\put(0,70){\framebox(50,40){T}}
\put(60,70){\oval(110,70)[b]}
\put(60,70){\oval(90,55)[b]}
\put(60,70){\oval(30,20)[b]}
\put(60,110){\oval(110,70)[t]}
\put(60,110){\oval(90,55)[t]}
\put(60,110){\oval(30,20)[t]}
\put(75,70){\line(0,1){40}}
\put(105,70){\line(0,1){40}}
\put(115,70){\line(0,1){40}}
\multiput(25,116)(6,0){3}{\makebox(0,0){.}}
\end{picture}
}
\caption{The closure link of a string link 
$T$\label{closure}}
\end{figure}
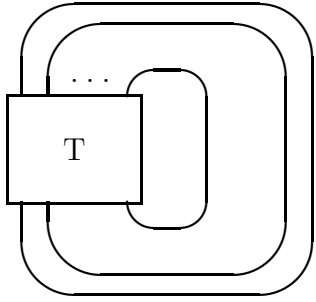

For a  string link $T$, one  defines the {\it closure link} of $T$ as
in the braid case, see Figure~\ref{closure}. 

Suppose that a {\it framed} string link $T$ is represented
by a tangle diagram on the plane, where the framing vector at every
boundary point is given by the vector perpendicular to the page and
 pointing to us. In the closure link of $T$, we suppose that the 
framing at every point outside $T$ is given by the vector 
perpendicular to the page and pointing to us. By this way we can
 identify the set of framings of a component of $T$ with $\Bbb Z$; say,
$C$ has framing $0$ if its closure has framing 0.
We say that $T$ is algebraically split if the closure link is
 algebraically split.

Let  $d_C(T)$ be obtained from $T$ by replacing component  $C$ by a
vertical line connecting the two boundary point of $C$ 
and  over-passing any other component of the tangle diagram. The
 new component is assumed to have the same framing (in $\Bbb Z$) as 
 $C$ has.

\begin{lem}
Suppose that a framed string link $T$  is
 {\em algebraically split}. Then $\hZ(T)-\hZ(d_C(T))
\in \cC(T)$ is {\rm i}-near $C$. 
\label{x13}
\end{lem}
\begin{pf}We suppose that $C$ is the last component of $T$, counting 
from left to right.
The other cases are similar (and also follow from this case and 
Lemma \ref{hZ}). Without loss of generality, we can assume that 
$C$ has framing 0.

Let $\hZ(T)=\exp(a+b)$, where $\varepsilon_C(a)=0$ 
and $\varepsilon_C(b)\not=0$. Both $a$ and $b$ are primitive.

Let $T'$ be obtained from $T$ by removing $C$. Then
(with $1\in\cP_1$)
$$\hZ(d_C(T))=\hZ(T')\otimes 1=\exp[\varepsilon_C(a+b)]\otimes 
1=\exp(\varepsilon_C(b))\otimes 1=\exp(b);$$
the first equality
follows from Lemma \ref{hZ}.

Hence  $\hZ(T)-\hZ(d_C(T))=\exp(a+b)-\exp(b)$ is a sum of terms,
each is a product with $a$ being one of the factor.

Since $\varepsilon_C(a)=0$, and $a$ is primitive, it follows that 
$a$ is a linear
combination of chord diagrams whose dashed graph is connected and 
has an 
external vertex on $C$. 

 Since the coefficients
of first degree chord diagrams 
in $\hZ(T)$ are the linking numbers or self-linking numbers (i.e.
framing numbers), we 
conclude that $a$ does not have any chord diagrams of degree 1.

It is clear that if $\xi$ is a chord diagram whose graph is 
connected, has an external
vertex on $C$ and $\xi$ has degree $\ge 2$, then $\xi$ has an internal
 vertex near $C$. Hence $a$ is i-near $C$.

It follows that  $\hZ(T)-\hZ(d_C(T))$ is {\rm i}-near $C$.
\end{pf}

\subsection{End of proof of Lemma \protect{\ref{lemma2}}}

An {\it $m$-marked Chinese character} is a vertex-oriented 
3-valent graph
whose external vertices are partitioned into $m$ sets $\Theta_1,\dots,
\Theta_m$. 
Consider the space $\cE_m$ spanned by $m$-marked Chinese
 characters, subject
to the AS and IHX relations.

Let $s$ be a subset of $\{1,2,\dots,m\}$.
An element $\xi\in\cE_m$ is i-{\it near} $s$ if $\xi$ is a 
linear
combination of $m$-marked Chinese characters, each has 
internal vertices
near $\Theta_j$ for every $j\in s$. Here an internal vertex 
is near $\Theta_j$
if it is connected to an external vertex in $\Theta_j$ by a 
dashed line.

Consider the linear mapping $\chi:\cE_m\to\cP_m$, 
defined as follows.
Suppose an $m$-marked Chinese character $\xi$ has $k_i$ 
external vertices in $\Theta_i,i=1,2,\dots,m$. There are 
$k_i!$ ways to put the vertices
from $\Theta_i$ on the $i$-th solid string; and each of the 
$k_1!\,k_2!\dots k_m!$ possibilities give us a chord 
diagram in $\cP_m$.
 Summing up all such
chord diagrams, we get $\chi(\xi)$.

It is known that $\chi$ is an isomorphism of vector
 spaces (see Theorem 9
in \cite{LM2}). The case
when $m=1$ was first proved by Kontsevich; for a 
detailed proof of the case
$m=1$ see 
 \cite{BN1} (Theorem 6).
The inverse of $\chi$ is obtained by a symmetrizing process. 

It follows from the definitions of $\chi$ and its inverse that

\begin{pro}For the $j$-th component $C$ of the support 
of $\cP_m$,
$\chi:\cE_m\to\cP_m$ maps isomorphically the subspace of \/ 
 {\rm i}-near $j$ elements to 
the subspace of \/ {\rm i}-near $C$ elements.\label{x12}
\end{pro}

An important fact is that the AS and IHX relations preserve 
the property to
have an internal vertex near $\Theta_j$. Hence the space
 spanned by {\it $m$-marked
Chinese characters} {\it  not} {\rm i}-near $j$ has 0 
intersection with the space of elements i-near $j$, for 
every $j\in\{1,2,\dots,m\}$.
 The similar
fact does not hold true for $\cP_m$, due to the STU relation. 
It follows that

\begin{pro}
Suppose that $\xi\in\cE_m$ is {\rm i}-near $j$ for each $j\in s$, 
where
$s$ is a subset of $\{1,2,\dots,m\}$, then $\xi$ is 
{\rm i}-near $s$.
\end{pro}

From this proposition and  Proposition \ref{x12} we have

\begin{pro}\label{x14}
If $\xi\in\cP_m$ is {\rm i}-near $C$ for each component
 $C$ of\/ $Y$, then
$\xi$ is {\rm i}-near $Y$.
\end{pro}

Now we prove Lemma \ref{lemma2}. We need to prove 
(\ref{equ:lemma2}). We will first prove that 
$\cZ(\tilde \delta_{L_2}(L))$
is i-near $L_2$, i.e.,  $\cZ(\tilde \delta_{L_2}(L))$
 is a linear combination of chord diagrams,
 each has internal vertices near every component of $L_2$. Since 
each internal vertex can be near at most 3 different solid components,
 and since $|L_2|=6n+1$, 
it then follows that $\cZ(\tilde \delta_{L_2}(L))$ has 
{\rm i}-filter $\ge 2n+1$. Hence
$\Omega_n[\cZ(\tilde \delta_{L_2}(L))]=0$.

 We can represent $L=(L_1,L_2)$ as 
the closure of an $m$-component string link $T$. 
We suppose $T$ is the union of $T_i$, $i=1,2$, 
where the components of $T_i$  
correspond to the 
components of $L_i$. For a sub-tangle $T_2'\subset T_2$,
 let $(T_1,T_2'')$ be the tangle 
obtained from $T=(T_1,T_2)$ by successively applying 
$d_C$ to $T$, where $C$ 
runs the set of all components of $T_2\setminus T_2'$ 
and $d_C$ is as in
 Lemma ~\ref{x13}.
Let $$\tilde \delta_{T_2}(T)=\sum_{T_2'\subset T_2}
(-1)^{ |T_2'|}(T_1,T_2'').$$

It is easy to see that for every component $C$ of 
$T_2$, the sum $\hZ
(\tilde \delta_{T_2}(T))$ can be represented as a
 finite sum of terms, each is of
 the form
$\hZ(V)-\hZ(d_C(V))$. Hence by
 Lemma \ref{x13},
$\hZ(\tilde \delta_{T_2}(T))\in\cP_m$ is {\rm i}-near 
$C$. Since this is
 true
for every component $C$ in $T_2$, Proposition \ref{x14} 
says that 
$\hZ(\tilde \delta_{T_2}(T))$ is i-near $T_2$.

Observe that $\tilde \delta_{T_2}(T)$ is a linear combination
of framed string links; and the closure of $\tilde \delta_{T_2}(T)$
is exactly $\tilde \delta_{L_2}(L)$.
By Lemma \ref{hZ} and property (\ref{multi}) we have
$$\cZ(\tilde \delta_{L_2}(L))=a[\hZ(\tilde 
\delta_{T_2}(T))\otimes 1]b,$$
where $1\in\cP_m$ is the chord diagram without 
dashed graph, and $a,b$ are some
linear combinations of chord diagrams. It follows that 
 $\cZ(\tilde \delta_{L_2}(L))$
is i-near $L_2$. As argued above, since $|L_2| \ge 6n+1$, we get
$\Omega_n[\cZ(\tilde \delta_{L_2}(L))]=0$. This completes
 the proof of 
(\ref{equ:lemma2}), and hence that of Lemma \ref{lemma2}.

\subsection{Proof of Lemma \protect{\ref{phu111}}}\label{x50}
 By Lemma \ref{hZ}, it suffices to consider the case when 
 $V_i$ does not
contain any vertical line. Then $\beta(V_i)=u-v$, where 
$v$ is 3 arcs not linked 
with each other, while $u$ is 3 arcs linked with each 
other like in the Borromean
link.

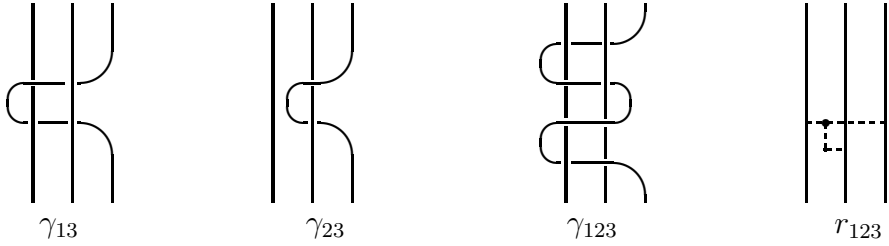
\begin{figure}[htp]
\centerline{
\begin{picture}(40,95)
\thicklines
\put(8,57.5){\oval(15,15)[l]}
\put(27,95){\oval(26,60)[br]}
\put(27,20){\oval(26,60)[tr]}
\put(10,20){\line(0,1){43.5}}
\put(10,95){\line(0,-1){28.5}}
\put(25,20){\line(0,1){75}}
\put(12,50){\line(1,0){11}}
\put(8,65){\line(1,0){14}}
\put(20,10){\makebox(0,0){$\gamma_{13}$}}
\end{picture}
\hskip 2cm
\begin{picture}(40,95)
\thicklines
\put(13,57.5){\oval(15,15)[l]}
\put(17,95){\oval(26,60)[br]}
\put(17,20){\oval(26,60)[tr]}
\put(15,20){\line(0,1){43.5}}
\put(15,95){\line(0,-1){28.5}}
\put(0,20){\line(0,1){75}}
\put(13,65){\line(1,0){4}}
\put(20,10){\makebox(0,0){$\gamma_{23}$}}
\end{picture}
\hskip 2cm
\begin{picture}(40,95)
\thicklines
\put(8,42.5){\oval(15,15)[l]}
\put(8,72.5){\oval(15,15)[l]}
\put(27,57.5){\oval(15,15)[r]}
\put(27,95){\oval(26,30)[br]}
\put(27,20){\oval(26,30)[tr]}
\put(10,20){\line(0,1){28.5}}
\put(10,52){\line(0,1){11.5}}
\put(10,95){\line(0,-1){28.5}}
\put(25,20){\line(0,1){13.5}}
\put(25,36.5){\line(0,1){11.5}}
\put(25,95){\line(0,-1){43.5}}
\put(12,35){\line(1,0){15}}
\put(8,50){\line(1,0){19}}
\put(8,65){\line(1,0){15}}
\put(12,80){\line(1,0){11.5}}
\put(20,10){\makebox(0,0){$\gamma_{123}$}}
\end{picture}
\hskip 2cm
\begin{picture}(30,95)
\thicklines
\put(0,20){\line(0,1){75}}
\put(15,20){\line(0,1){75}}
\put(30,20){\line(0,1){75}}
\put(0,50){\dashbox{2}(30,0){}}
\put(7.5,40){\dashbox{2}(7.5,0){}}
\put(7.5,40){\dashbox{2}(0,10){}}
\put(7.5,50){\circle*{2.5}}
\put(20,10){\makebox(0,0){$r_{123}$}}
\end{picture}
}
\caption{$\gamma_{13},\gamma_{23}$, 
$\gamma_{123}=[\gamma_{13}^{-1},\gamma_{23}]$, and $r_{123}$
\label{gamma}}
\end{figure}

These two tangles $u$ and $v$ can also be described as follows.
 The tangle  $u$ is obtained from the tangle
in Figure \ref{x24} by replacing the box by $\gamma_{123}$, 
where $\gamma_{123}$ is
the tangle described in Figure \ref{gamma}; and $v$ is 
obtained by replacing the
box by the trivial tangle on 3 strands. 
The orientations of $u,v$ are
chosen so that all the lines in the box are  upward. 
\begin{figure}[htp]
\centerline{
\begin{picture}(90,130)
\thicklines
\put(45,100){\oval(90,60)[t]}
\put(45,100){\oval(70,42)[t]}
\put(45,100){\oval(50,25)[t]}
\put(-10,80){\framebox(40,20){box}}
\put(90,100){\line(0,-1){20}}
\put(80,100){\line(0,-1){28}}
\put(70,100){\line(0,-1){38}}
\put(70,58){\line(0,-1){36}}
\put(80,68){\line(0,-1){38}}
\put(90,80){\line(-1,-1){80}}
\put(0,80){\line(0,-1){80}}
\put(10,80){\line(0,-1){50}}
\put(20,80){\line(0,-1){20}}
\put(40,0){\line(-1,1){13}}
\put(80,0){\line(-1,1){13}}
\put(90,0){\line(-1,1){18}}
\put(50,0){\line(1,1){30}}
\put(10,30){\line(1,-1){13}}
\put(20,60){\line(1,-1){23}}
\put(47,33){\line(1,-1){16}}
\end{picture}
}
\caption{}
\label{x24}
\end{figure}

Let $\gamma_{13},\gamma_{23}$ be elements of the pure
 braid group $P_3$ on 3 strands
depicted in Figure \ref{gamma}. We suppose the orientation of 
each strand is upwards.
An important observation is that $\gamma_{123}=
\gamma_{13}^{-1}\gamma_{23}
\gamma_{13}\gamma_{23}^{-1}$, a commutator.

An easy exercise in the theory of the Kontsevich integral is
that, for ${ij}=12$ or $13$,
\begin{equation}\hZ(\gamma_{ij})=\exp(r_{ij}+\xi_{ij}),
\label{x20}\end{equation}
where  $r_{ij}$ is the chord diagram of
 degree 1 with one dashed line
connecting the $i$-th and $j$-th strands, and $\xi_{ij}$ 
is of degree $\ge 2$.
Hence for $\gamma_{123}=\gamma_{13}^{-1}\gamma_{23}
\gamma_{13}\gamma_{23}^{-1}$, one has
$$\hZ(\gamma_{123})=1-r_{13}r_{23}+r_{23}r_{13}+
(\text{elements of degree} \ge 3),$$
By the STU relation $-r_{13}r_{23}+r_{23}r_{13}=
r_{123}$, where $r_{123}$ is the chord
diagram in Figure \ref{gamma}. Hence we have
$$\hZ(\gamma_{123})=\exp(r_{123}+\xi),$$
where $\xi$ is primitive and  of degree $\ge 3$. It follows
from Proposition \ref{ifilter} that
\begin{equation}
\hZ(\gamma_{123})=1+r_{123}+(\text{elements of 
{\rm i}-filter $\ge 2$}).
\label{x21}
\end{equation}
Note that if we replace the box in Figure \ref{x24}
 by $r_{123}$, then we get $\xi(V_i)$.
Hence 
$$\hZ(\beta(V_i))=\hZ(u)-\hZ(v)=\xi(V_i) + 
\text{(elements of {\rm i}-filter $\ge 2$)}.$$
This completes the proof.


\begin{thebibliography}{99999}
\baselineskip15pt

\bibitem[AS1]{AS1}
S. Axelrod and I. M. Singer,
{\it Chern-Simons perturbation theory},
Proc. XXth DGM Conference (New York, 1991)
World Scientific, 3--45 (1992).

\bibitem[AS2]{AS2}
\bysame,
{\it Perturbative aspects of Chern-Simons topological quantum
field theory II},
J. Diff. Geom. {\bf 39}, 173--213 (1994).


\bibitem[BN1]{BN1} D. Bar-Natan,
{\it On the Vassiliev knot invariants},
Topology, {\bf 34}, 423--472 (1995).

  \bibitem[BN2]{BN2} 
 \bysame,
{\it Non-Associative tangles},
Harvard University preprint, 1993.

\bibitem[BL]{BL} J.S. Birman
and X.S. Lin, {\it Knot polynomials
 and Vassiliev's invariants}, { Invent. Math.},
{\bf 111} (1993), pp.
 225--270.

\bibitem[GL1]{GL}
S. Garoufalidis and J. Levine,
{\it On finite type 3-manifold invariants II}, preprint 1995.

\bibitem[GL2]{GL2}
\bysame,
{\it On finite type 3-manifold invariants IV:
 comparison of definitions}, preprint.

\bibitem[GO1]{GO}
S. Garoufalidis and T. Ohtsuki,
{\it On finite type 3-manifold invariants III: 
manifold weight systems}, 
preprint, August 1995.

\bibitem[GO2]{GO2}
\bysame,
{\it On finite type 3-manifold invariants V}, preprint, 1995.


\bibitem[Gus]{Gus} M. N. Gusarov,
{\it On $n$-equivalence of knots and invariants of finite degrees}, 
in 
``topology of manifolds and varieties", edited by O. Viro, Adv.
 Soviet Math.,
 {\bf 18}, 1994.

\bibitem[Hab]{Hab}
N. Habegger, {\it Finite type 3-manifold invariants: a proof 
of a conjecture of
Garoufalidis}, preprint, 1995.

\bibitem[Hos]{Hos} J. Hoste,
{\it A formula for Casson's invariant},
Trans. Amer. Math. Soc., {\bf 297} (1986), 547--562.


\bibitem[Kas]{Kas}C. Kassel, {\it Quantum groups},
  Graduate Texts in Mathematics {\bf 155},
Springer-Verlag, in 1994.

\bibitem[Koh]{Kohno} T. Kohno,
{\it Vassiliev invariants and de-Rham complex on the space of knots},
preprint, Tokyo university, 1994.


\bibitem[Ko1]{Ko1}
M. Kontsevich,
{\it Vassiliev's knot invariants},
Advances in Soviet Mathematics, {\bf 16}, 137--150 (1993).

\bibitem[Ko2]{Ko2}
\bysame,
{\it Feynman diagrams and low-dimensional topology},
preprint 1992.

\bibitem[LM1]{LM1}
T. T. Q. Le and J. Murakami,
{\it Representations of the category of tangles by 
Kontsevich's iterated
integral},
Commun. Math. Phys., {\bf 168}, 535--562 (1995).

\bibitem[LM2]{LM2} 
\bysame,
{\it The universal Vassiliev-Kontsevich invariant for 
framed oriented links}, to appear in Compositio Math.

\bibitem[LM3]{LM3} 
\bysame,
{\it Parallel version of the universal Vassiliev-Kontsevich 
invariant},
preprint 1994.

\bibitem[LMO]{LMO}
T. T. Q. Le, J. Murakami and T. Ohtsuki,
{\it On a universal quantum invariant of 3-manifolds},
preprint q-alg/9512002, 1995.

\bibitem[LMMO]{LMMO}
T. T. Q. Le, H. Murakami, J. Murakami and T. Ohtsuki,
{\it A three-manifold invariant derived from the universal
Vassiliev-Kontsevich invariant},
Proc. Japan Acad., {\bf 71}, Ser. A, 125--127 (1995).

\bibitem[Lin]{Lin} X. S. Lin,
 {\it Finite type invariants of integral homology 3-sphere:
 a survey}, preprint, October 1995.

\bibitem[Ng]{Ng} Ka Yi Ng,
{\it Groups of ribbon knots}, preprint, Columbia university, 1995.

\bibitem[NgS]{SN} Ka Yi Ng, T. Stanford,
{\it On Gusarov's groups of knots}, preprint 1995.

\bibitem[Oh1]{Ohtsuki}
T. Ohtsuki,
{\it Finite type invariant of integral homology 3-spheres},
 to appear in
Jour. of Knot theory Rami..

\bibitem[Oh2]{Oht2}
\bysame,
{\it A polynomial invariant of rational homology 3-spheres},
to appear in Invent. Math.
\bibitem[Oh3]{Oht3}\bysame,
{\it Invariants of 3-manifolds derived from universal 
invariants of framed links},
Math. Proc. Camb. Phil. Soc., {\bf 117} (1995), 259--273.

\bibitem[RT]{RT2}
N. Reshetikhin and V. G. Turaev,
{\it Invariants of $3$-manifolds via link polynomials
 and quantum groups},
Invent.\ Math.  {\bf 103}, 547--597 (1991).

\bibitem[Roz]{Roz}
L. Rozansky,
{\it The trivial connection contribution to Witten's invariant and
finite type invariants of rational homology spheres},
preprint, 1995.
\bibitem[Ser]{Ser} J. P. Serre,{ \it  Lie algebras 
and Lie groups,}, W.
 A. Benjamin Inc., New York, Amsterdam, 1965.

 \bibitem[Sta]{Sta} T. Stanford, 
{\it Braid commutators and Vassiliev invariants}, 
to appear in Pacific J. Math..

\bibitem[Tur]{Tur} V.~G. Turaev, {\it Quantum invariants of 
knots and 3-manifolds}, de Gruyter Studies in Mathematics 18, 
Walter de Gruyter, Berlin New York 1994.

 \bibitem[Vas]{Vas} 
V.~A.~Vassiliev,
{\it Cohomology of knot spaces},
 in ``Theory of Singularities and Its Applications'' 
(ed. V.I. Arnold),
 Advances in Soviet Mathematics, AMS {\bf 1} (1990),  23--69.

\bibitem[Wit]{Witten}
E. Witten,
{\it Quantum field theory and the Jones polynomial},
Commun. Math. Phys. {\bf 121}, 360--379 (1989).

\end{thebibliography}
        \end{document}